\documentclass[11pt]{article}
\usepackage{amssymb}
\usepackage{amsfonts}
\usepackage{amsmath}
\usepackage{graphicx}

\newtheorem{theorem}{Theorem}

\newtheorem{lemma}[theorem]{Lemma}

\newtheorem{proposition}[theorem]{Proposition}
\newtheorem{remark}[theorem]{Remark}

\newenvironment{proof}[1][Proof]{\noindent\textbf{#1.} }{\ \rule{0.5em}{0.5em}}

\makeatletter\@addtoreset {equation}{section}\makeatother

\begin{document}

\title{\bf Incompressible viscous fluid flows \\ in a thin spherical shell}
\author{Ranis N. Ibragimov and Dmitry E. Pelinovsky \\
{\small Department of Mathematics, McMaster University, Hamilton,
Ontario, Canada, L8S 4K1} }

\maketitle

\begin{abstract}
Linearized stability of incompressible viscous fluid flows in a
thin spherical shell is studied by using the two-dimensional
Navier--Stokes equations on a sphere. The stationary flow on the
sphere has two singularities (a sink and a source) at the North
and South poles of the sphere. We prove analytically for the
linearized Navier--Stokes equations that the stationary flow is
asymptotically stable. When the spherical layer is truncated
between two symmetrical rings, we study eigenvalues of the
linearized equations numerically by using power series solutions
and show that the stationary flow remains asymptotically stable
for all Reynolds numbers.
\end{abstract}

\section{Introduction}

The Navier-Stokes (NS) equations for an incompressible viscous
fluid are the fundamental governing equations of fluid mechanics.
In many cases, exact solutions can be constructed to these
equations \cite{DR06} and spectral and nonlinear stability of
these exact solutions can be analyzed \cite{DR81}. Our work
addresses stability of exact solutions for the NS equations in
spherical coordinates.

The three-dimensional NS equations in a thin rotating spherical
shell describe large-scale atmospheric dynamics that plays an
important role in the global climate control and weather
prediction \cite{LTW92,LTW92b} (see also review in \cite{Gill}).
It was rigorously proved by Temam \& Ziane \cite{TZ97} that the
average of the longitudinal velocity in the radial direction
converges to the strong solution of the two-dimensional NS
equation on a sphere as the thickness of the spherical shell goes
to zero. The latter model has been used in geophysical fluid
dynamics since middle of the last century \cite{Russian_Text}.

The treatment of the geometric singularity in spherical coordinates
has for many years been a difficulty in the development of numerical
simulations for oceanic and atmospheric flows around the Earth.
Blinova \cite{Blinova1,Blinova2} represented solutions in the
inviscous case by the eigenfunction expansions in spherical
harmonics. Vorticity equations were considered by Ben-Yu with the
spectral method \cite{GB}. More recent work of Furnier et al.
\cite{FB} applied the spectral-element method to the axis-symmetric
solutions (see \cite{JB,MC,W} for other applications of the spectral
methods in spherical coordinates). Finally, point vortex motion on a
sphere was modeled by ordinary differential equations for vortex
centers in Boatto \& Cabral \cite{B06} and Crowdy \cite{Crowdy}.

We address the three-dimensional NS equations for an incompressible
viscous fluid,
\begin{equation}
\label{1.1} \left\{ \begin{array}{lll} & \frac{\partial
\mathbf{u}}{\partial t}+\left( \mathbf{u\cdot \nabla }\right)
\mathbf{u}-\nu \Delta \mathbf{u}+\nabla p=0, \quad & {\bf x} \in
\Omega, \; t \in \mathbb{R}_+, \\ & \mathbf{\nabla }\cdot
\mathbf{u}=0, \quad & {\bf x} \in \Omega, \; t \in \mathbb{R}_+, \\
& {\bf u} |_{t = 0} = {\bf u}_0, \quad & {\bf x} \in \Omega,
\end{array} \right.
\end{equation}
in a thin spherical shell $\Omega = \{ {\bf x} \in \mathbb{R}^3 :
1 < |{\bf x}| < 1 + \varepsilon \}$ with $\varepsilon \to 0$,
subject to the boundary conditions
\begin{equation}
\label{1.3} {\bf u} \cdot {\bf n} = 0, \qquad \nabla {\bf u} \times
{\bf n} = {\bf 0}, \qquad {\bf x} \in \partial \Omega.
\end{equation}
Here $\mathbf{u} : \Omega \times \mathbb{R}_+ \mapsto
\mathbb{R}^3$ is the velocity vector, $p : \Omega \times
\mathbb{R}_+ \mapsto \mathbb{R}$ is the ratio of the pressure to
constant density, $\nu$ is the kinematic viscosity, ${\bf n}$ is
the normal vector to the boundary $\partial \Omega$ of the
spherical shell $\Omega$ and $\mathbf{u}_{0} : \Omega \mapsto
\mathbb{R}^3$ is a given initial condition. Although Coriolis and
gravity forces may be dynamically significant in oceanographic
applications, our model is considered in a non-rotating reference
frame and without external forces. The effects of rotation and
gravity can be included into the model but they do not
substantially alter the physical picture that emerges from the NS
equations (\ref{1.1}).

We employ the spherical coordinates $(r,\theta,\phi)$ with the
velocity vector ${\bf u} = u_r {\bf e}_r + u_{\theta} {\bf
e}_{\theta} + u_{\phi} {\bf e}_{\phi}$, where $({\bf e}_r,{\bf
e}_{\theta},{\bf e}_{\phi})$ are basic orthonormal vectors along
the spherical coordinates. For completeness, we reproduce the
three-dimensional NS equations (\ref{1.1}) in spherical
coordinates \cite{Bachelor}:
\begin{eqnarray*}
&& \frac{\partial u_r}{\partial t} + u_r \frac{\partial
u_r}{\partial r} + \frac{u_{\theta}}{r} \frac{\partial
u_r}{\partial \theta} + \frac{u_{\phi}}{r \sin \theta}
\frac{\partial u_r}{\partial \phi} - \frac{u_{\theta}^2 +
u_{\phi}^2}{r} = - \frac{\partial p}{\partial r} + \nu \left(
\Delta u_r + \frac{2}{r} \frac{\partial
u_r}{\partial r} + \frac{2 u_r}{r^2}  \right), \\
\nonumber && \frac{\partial u_{\theta}}{\partial t} + u_r
\frac{\partial u_{\theta}}{\partial r} + \frac{u_{\theta}}{r}
\frac{\partial u_{\theta}}{\partial \theta} + \frac{u_{\phi}}{r
\sin \theta} \frac{\partial u_{\theta}}{\partial \phi} + \frac{u_r
u_{\theta}}{r} - \frac{u_{\phi}^2 \cot \theta}{r} = - \frac{1}{r}
\frac{\partial p}{\partial \theta} \\
&& \phantom{texttexttexttexttexttext} + \nu \left( \Delta
u_{\theta } + \frac{2}{r^2} \frac{\partial u_r}{\partial \theta} -
\frac{u_{\theta }}{r^2 \sin^2 \theta} - \frac{2\cos \theta}{r^2
\sin ^2 \theta } \frac{\partial u_{\phi }}{\partial \phi} \right),
 \\ \nonumber && \frac{\partial u_{\phi}}{\partial t} +
u_r \frac{\partial u_{\phi}}{\partial r} + \frac{u_{\theta}}{r}
\frac{\partial u_{\phi}}{\partial \theta} + \frac{u_{\phi}}{r \sin
\theta} \frac{\partial u_{\phi}}{\partial \phi} + \frac{u_r
u_{\phi}}{r} + \frac{u_{\theta} u_{\phi} \cot \theta}{r} = -
\frac{1}{r \sin \theta} \frac{\partial p}{\partial \phi} \\
&& \phantom{texttexttexttexttexttext} + \nu \left( \Delta u_{\phi
} + \frac{2}{r^2 \sin \theta} \frac{\partial u_r}{\partial \phi} +
\frac{2\cos \theta }{r^2 \sin^2 \theta } \frac{\partial u_{\theta
}}{\partial \phi }-\frac{u_{\phi }}{r^2 \sin^{2}\theta }\right), \\
&& \frac{1}{r^2} \frac{\partial}{\partial r} \left( r^2 u_r
\right) + \frac{1}{r \sin \theta} \frac{\partial }{\partial \theta
}\left( \sin \theta u_{\theta} \right) + \frac{1}{r \sin \theta}
\frac{\partial u_{\phi }}{\partial \phi } = 0,
\end{eqnarray*}
where
\begin{eqnarray*}
\Delta = \frac{1}{r^2} \frac{\partial }{\partial r} \left( r^2
\frac{\partial}{\partial r} \right) + \frac{1}{r^2 \sin \theta}
\frac{\partial}{\partial \theta} \left( \sin \theta
\frac{\partial}{\partial \theta} \right) + \frac{1}{r^2 \sin^2
\theta} \frac{\partial^2}{\partial \phi^2}
\end{eqnarray*}
is the Laplacian in spherical coordinates and the initial and
boundary conditions are not written. One can check by direct
differentiation that there exists an exact stationary solution to
the three-dimensional NS equations in spherical coordinates:
\begin{equation}
\label{exact-solution} u_r = 0, \quad u_{\theta} = \frac{\alpha}{r
\sin \theta}, \quad u_{\phi} = 0, \quad p = \beta -
\frac{\alpha^2}{2 r^2 \sin^2 \theta},
\end{equation}
where $(\alpha,\beta)$ are arbitrary parameters. The stationary
solution (\ref{exact-solution}) describes fluid motion tangential
to a sphere of any given radius $r$. The stationary flow has two
pole singularities at $\theta = 0$ and $\theta = \pi$. The
singularities correspond to the source and sink of the velocity
vector at the North and South poles of the spherical shell
$\Omega$: the fluid is injected at the North pole from an external
source and it leaks out at the South pole to an external sink.

In the limit $\varepsilon \to 0$, the non-stationary
three-dimensional fluid flow is confined on a sphere $S$ of unit
radius parameterized by the polar (latitude) angle $\theta$ and
azimuthal (longitude) angle $\phi$,
\begin{equation}
\label{domain-complete} S =\left\{ \left( \theta ,\phi \right)
,\text{ \ }0\leqslant \theta \leqslant \pi ,\text{ }0\leqslant \phi
<2\pi \right\}.
\end{equation}
Since the velocity vector ${\bf u}$ and the pressure $p$ in the NS
equations (\ref{1.1}) are coupled together by the
incompressibility constraint $\nabla \cdot {\bf u} = 0$, it is
difficult to analyze the full set of three-dimensional equations.
A common approach to simplify the problem is to use the artificial
methods such as the pressure stabilization and projections
\cite{Shen}. The error estimate of the pressure stabilization and
projection methods is not however mathematically precise. Instead,
we shall use the result of the Theorem B in \cite{TZ97}, which
states that provided the function ${\bf u}_0(r,\theta,\phi)$ is
smooth enough, the strong global solution ${\bf
u}(r,\theta,\phi,t)$ of the three-dimensional NS equations
converges as $\varepsilon \to 0$ to the strong unique global
solution ${\bf v}(\theta,\phi,t)$ of the two-dimensional NS
equations on the sphere, where
$$
{\bf v}(\theta,\phi,t) = \lim\limits_{\varepsilon \to 0}
\frac{1}{\varepsilon} \int_1^{1+\varepsilon} r {\bf
u}(r,\theta,\phi,t) dr = (0,v_{\theta},v_{\phi}).
$$
The vector ${\bf v}(\theta,\phi,t)$ is interpreted as the average
velocity with respect to the radial coordinate $r$. The
two-dimensional NS equations on a sphere $S$ in spherical angles
$(\theta,\phi)$ are written explicitly as follows \cite{TZ97}:
\begin{eqnarray*}
&& \frac{\partial v_{\theta}}{\partial t} + v_{\theta}
\frac{\partial v_{\theta}}{\partial \theta} + \frac{v_{\phi}}{\sin
\theta} \frac{\partial v_{\theta}}{\partial \phi} - v_{\phi}^2
\cot \theta = - \frac{\partial p}{\partial \theta} + \nu \left(
\Delta_S v_{\theta } - \frac{v_{\theta }}{\sin^2 \theta} -
\frac{2\cos \theta}{\sin ^2 \theta } \frac{\partial v_{\phi
}}{\partial \phi} \right),  \\ && \frac{\partial
v_{\phi}}{\partial t} + v_{\theta} \frac{\partial
v_{\phi}}{\partial \theta} + \frac{v_{\phi}}{\sin \theta}
\frac{\partial v_{\phi}}{\partial \phi} + v_{\theta} v_{\phi} \cot
\theta = - \frac{1}{\sin \theta} \frac{\partial p}{\partial \phi}
+ \nu \left( \Delta_S v_{\phi } + \frac{2\cos \theta }{\sin^2
\theta } \frac{\partial v_{\theta
}}{\partial \phi} - \frac{v_{\phi }}{\sin^{2}\theta }\right), \\
&& \frac{1}{\sin \theta} \frac{\partial }{\partial \theta }\left(
\sin \theta v_{\theta} \right) + \frac{1}{\sin \theta}
\frac{\partial v_{\phi }}{\partial \phi } = 0,
\end{eqnarray*}
where $\Delta_S$ is the Laplace-Beltrami operator in spherical
angles
\begin{eqnarray*}
\Delta_S =\frac{1}{\sin \theta }\frac{\partial }{\partial \theta
}\left( \sin \theta \frac{\partial }{\partial \theta }\right) +\frac{1}{\sin ^{2}\theta }%
\frac{\partial ^{2}}{\partial \phi ^{2}}.
\end{eqnarray*}
Note that no boundary conditions are specified for the vector
${\bf v}(\theta,\phi,t)$ on sphere $S$, while the initial
condition ${\bf v}|_{t = 0} = {\bf v}_0$ on $S$ is not written.
For the purposes of our work, we rewrite the two-dimensional NS
equations on the sphere $S$ in an equivalent form:
\begin{eqnarray}
&& \frac{\partial v_{\theta }}{\partial t} - \frac{v_{\phi}
\omega}{\sin \theta} + \frac{\partial q}{\partial \theta } = \nu
\left( \Delta_S v_{\theta} - \frac{v_{\theta}}{\sin^2\theta }
-\frac{2\cos \theta }{\sin ^{2}\theta} \frac{\partial v_{\phi }}{\partial \phi} \right),
\label{2.1} \\
&& \frac{\partial v_{\phi}}{\partial t} + \frac{v_{\theta}
\omega}{\sin \theta} + \frac{1}{\sin \theta} \frac{\partial
q}{\partial \phi } = \nu \left( \Delta_S v_{\phi }+\frac{2\cos
\theta }{\sin ^{2}\theta } \frac{\partial v_{\theta
}}{\partial \phi }-\frac{v_{\phi }}{\sin^{2}\theta }\right),  \label{2.2} \\
&& \frac{\partial }{\partial \theta }\left( \sin \theta v_{\theta}
\right) + \frac{\partial v_{\phi }}{\partial \phi } = 0, \label{2.3}
\end{eqnarray}
where $q$ is a static (stagnation) pressure and $\omega$ is the
vorticity:
\begin{equation}
\label{vorticity-renorm} q = p + \frac{1}{2} \left( v_{\theta}^2 +
v_{\phi}^2 \right), \qquad \omega = \frac{\partial }{\partial
\theta }\left( \sin \theta \; v_{\phi } \right) -\frac{\partial
v_{\theta }}{\partial \phi}.
\end{equation}
The stationary solution (\ref{exact-solution}) corresponds to the
exact stationary solution of the two-dimensional NS equations
(\ref{2.1})--(\ref{2.3}) on the unit sphere $S$:
\begin{equation}
\label{exact-solution-2D} v_{\theta} = \frac{\alpha}{\sin \theta},
\qquad v_{\phi} = 0, \qquad q = \beta,
\end{equation}
where $(\alpha,\beta)$ are arbitrary parameters.

We shall also consider the situation where the external source and
sink singularities at $\theta = 0$ and $\theta = \pi$ are excluded
from the domain of the NS equations (\ref{2.1})--(\ref{2.3}). For
instance, we shall consider the truncated domain in the form of the
spherical layer
\begin{equation}
\label{domain-truncated} S_0 = \left\{ (\theta,\phi) : \quad
\theta_0 \leq \theta \leq \pi - \theta_0, \;\; 0 \leq \phi \leq
2\pi \right\},
\end{equation}
where $0 < \theta_0 < \frac{\pi}{2}$. Without loss of generality,
the spherical layer $S_0$ is truncated symmetrically at the two
rings located in the Northern and Southern semi-spheres such that
the stationary flow (\ref{exact-solution-2D}) is free of pole
singularities in $S_0$. In other words, without dipping into
details on how the fluid flow is injected on the sphere and is
collected from the sphere in a neighborhood of the North and South
poles, we will study how the fluid leaks from the Northern
semi-sphere to the Southern semi-sphere along the spherical layer
(\ref{domain-truncated}). In this context, the stationary solution
(\ref{exact-solution-2D}) is interpreted as the mass conservation
law which is obtained by integrating the free divergence condition
(\ref{2.3}).

We are interested in spectral stability of the stationary fluid flow
(\ref{exact-solution-2D}). In the case of $S$ when the singularities
are included, we prove analytically that the linearized NS equations
(\ref{2.1})--(\ref{2.3}) about the stationary solution
(\ref{exact-solution-2D}) are asymptotically stable. In the case of
$S_0$ when the singularities are excluded, the asymptotical
stability of the stationary flow can only be proved for the case
$\nu = \infty$, that is in the limit of zero Reynolds numbers. By
using the power series expansions, we approximate solutions
numerically and show that the stationary flow remains asymptotically
stable for all Reynolds numbers.

Our paper is structured as follows. Section 2 introduces the
linearization of the two-dimensional NS equations
(\ref{2.1})--(\ref{2.3}) at the stationary solution
(\ref{exact-solution-2D}) and discusses boundary conditions for the
perturbation vector. Analytical results on location of the spectrum
of the linearized problem are reported in Section 3 for
symmetry-breaking ($\phi$-dependent) perturbations and in Section 4
for symmetry-preserving ($\phi$-independent) perturbations.
Numerical results on computations of eigenvalues of the linearized
problem are described in Section 5 for symmetry-breaking
perturbations and in Section 6 for symmetry-preserving
perturbations. Section 7 discusses applications.

\section{Linearized equations and separation of variables}

Without loss of generality, we consider the stationary solution
(\ref{exact-solution-2D}) with $\alpha = 1$ and $\beta = 0$. The
presence of arbitrary parameters ($\alpha,\beta$) introduces
time-independent (neutral) modes of the linearized equations,
which we will also account for in this section. We consider an
infinitesimal time-dependent perturbations of the stationary flow
with $\alpha = 1$ and $\beta = 0$ in the form
\begin{equation}
v_{\theta} = \frac{1}{\sin \theta} + U(\theta,\phi) e^{\lambda t},
\qquad v_{\phi} = V(\theta,\phi) e^{\lambda t}, \qquad q =
Q(\theta,\phi) e^{\lambda t}, \label{2.6}
\end{equation}
where $\lambda \in \mathbb{C}$ is a parameter, such that
perturbations with ${\rm Re}(\lambda) > 0$ imply spectral
instability of the stationary flow. If ${\rm Re}(\lambda) < 0$ for
all perturbations, the stationary flow is asymptotically stable,
while if ${\rm Re}(\lambda) = 0$ for some perturbations and ${\rm
Re}(\lambda) < 0$ for all other perturbations, the stationary flow
is stable in the sense of Lyapunov.

By neglecting the quadratic terms of the perturbation, we linearize
the NS equations (\ref{2.1})--(\ref{2.3}) with the expansion
(\ref{2.6}) to the form:
\begin{eqnarray}
&& \lambda U + \frac{\partial Q}{\partial \theta } = \nu \left(
\Delta_S U - \frac{U}{\sin ^{2}\theta} - \frac{2\cos \theta }{\sin^2
\theta } \frac{\partial V}{\partial \phi }\right), \label{2.7} \\
&& \lambda V + \frac{1}{\sin^2 \theta} \left(
\frac{\partial}{\partial \theta} \left( \sin \theta \; V \right) -
\frac{\partial U}{\partial \phi} \right) + \frac{1}{\sin \theta}
\frac{\partial Q}{\partial \phi } = \nu \left( \Delta_S V +
\frac{2\cos \theta }{\sin2 \theta }\frac{\partial U}{\partial \phi}
-\frac{V}{\sin ^{2}\theta }\right) , \label{2.8} \\
&& \frac{\partial}{\partial \theta} \left( \sin \theta \; U \right)
+ \frac{\partial V}{\partial \phi }=0.  \label{2.9}
\end{eqnarray}
Perturbation terms of the velocity vector must satisfy some boundary
conditions in the domains $S_0$ or $S$. It is naturally to assume
that the velocity vector is periodic with respect to the angle
$\phi$:
\begin{equation}
\label{bc-periodic} U(\theta,\phi+2\pi) = U(\theta,\phi), \quad
V(\theta,\phi+2\pi) = V(\theta,\phi).
\end{equation}
Therefore, we look for Fourier series solutions of the system
(\ref{2.7})--(\ref{2.9}):
\begin{equation}
U(\theta,\phi) = \sum_{k \in \mathbb{Z}} U_k(\theta) e^{i k \phi},
\quad V(\theta,\phi) = \sum_{k \in \mathbb{Z}} V_k(\theta) e^{i k
\phi}, \quad Q(\theta,\phi) = \sum_{k \in \mathbb{Z}} Q_k(\theta)
e^{i k \phi}.
\end{equation}
We also require that the components $(U_k,V_k)$ of the velocity
vector be square integrable in $S_0$ or $S$ with respect to the
spherical weight:
\begin{equation}
\label{square-integrable} \int_{\theta_0}^{\pi-\theta_0} \left(
|U_k|^2 + |V_k|^2 \right) \sin \theta d \theta < \infty,
\end{equation}
where $0 \leq \theta_0 < \pi/2$. When the domain is the truncated
spherical shell $S_0$, we require that components the velocity
vector vanish at the regular end points of the domain:
\begin{equation}
\label{bc-dirichlet-0} U_k(\theta_0) = U_k(\pi - \theta_0) =
V_k(\theta_0) = V(\pi-\theta_0) = 0.
\end{equation}
The complete sphere $S$ with the singular end points will be
considered in the limit $\theta_0 \to 0$. We require that the
components of the vorticity in (\ref{vorticity-renorm}) vanish at
the singular end points of the domain:
\begin{equation}
\label{bc-dirichlet} \lim_{\theta \to 0} U_k(\theta) =
\lim_{\theta \to \pi} U_k(\theta) = \lim_{\theta \to 0} \sin
\theta V_k(\theta) = \lim_{\theta \to \pi} \sin \theta V_k(\theta)
= 0.
\end{equation}
It will be clear later that separation of variables is different
between the cases $k = 0$ and $k \neq 0$. We say that the
correction terms with $k = 0$ represent {\em symmetry-preserving}
perturbations of the stationary flow (\ref{2.6}), while the
correction terms with $k \neq 0$ represent {\em symmetry-breaking}
perturbations.

{\bf Case $k \neq 0$:} It follows from the divergence-free condition
(\ref{2.9}) that one can introduce the stream function
$\Psi_k(\theta)$ for the velocity vector $(U_k,V_k)$ as follows:
\begin{equation}
U_k = \frac{i k}{\sin \theta} \Psi_k(\theta), \qquad V_k =
-\Psi'_k(\theta). \label{stream-function-representation}
\end{equation}
The system of linearized equations (\ref{2.7})--(\ref{2.8}) reduces
to the coupled ODE system  for $\Psi_k(\theta)$ and $Q_k = i k
P_k(\theta)$:
\begin{eqnarray}
\frac{d}{d\theta} P_k & = & \frac{1}{\sin \theta} \left( \nu
\Delta_k \Psi_k - \lambda \Psi_k \right), \label{system-1} \\
\frac{k^2}{\sin \theta} P_k & = & \frac{d}{d\theta}\left( \nu
\Delta_k \Psi_k - \lambda \Psi_k \right) - \frac{1}{\sin \theta}
\Delta_k \Psi_k,\label{system-2}
\end{eqnarray}
where
\begin{equation}
\label{Laplacian-k} \Delta_k = \frac{d^2}{d \theta^2} + \frac{\cos
\theta}{\sin \theta} \frac{d}{d\theta} - \frac{k^2}{\sin^2 \theta}.
\end{equation}
Let $\Phi_k = \Delta_k \Psi_k$ be a new variable. Then, the
variable $P_k$ can be excluded from the system
(\ref{system-1})--(\ref{system-2}), such that the system reduces
to a closed second-order ODE:
\begin{equation}
\label{eigenvalue-equation} \nu \Delta_k \Phi_k -
\frac{\Phi_k'}{\sin \theta} = \lambda \Phi_k.
\end{equation}
Besides the relations (\ref{system-1}) and (\ref{system-2}) between
the pressure $P_k$, the stream function $\Psi_k$ and the vorticity
$\Phi_k$, we note another relation between these components:
\begin{equation}
\label{second-equation} \Phi_k = \Delta_k \Psi_k = \sin^2 \theta
\Delta_k P_k.
\end{equation}
Due to the boundary conditions (\ref{bc-dirichlet-0}) and the
representation (\ref{stream-function-representation}), the
solution $\Psi_k(\theta)$ for the truncated spherical layer $S_0$
is defined on a closed interval $\theta_0 \leq \theta \leq \pi -
\theta_0$ for $0 < \theta_0 < \pi/2$ subject to the boundary
conditions
\begin{equation}
\Psi_k(\theta_0) = \Psi_k'(\theta_0) = \Psi_k(\pi-\theta_0) =
\Psi_k'(\pi-\theta_0) = 0. \label{X0}
\end{equation}
Since $\theta = 0$ and $\theta = \pi$ are singular points of the
system (\ref{eigenvalue-equation})--(\ref{second-equation}) when
$\theta_0 \to 0$, the solution $\Psi_k(\theta)$ for the complete
sphere $S$ is defined on an open interval $0 < \theta < \pi$
satisfying the boundary conditions from (\ref{bc-dirichlet}) and
(\ref{stream-function-representation}):
\begin{equation}
\lim_{\theta \to 0} \Psi_k(\theta) = \lim_{\theta \to 0} \sin \theta
\Psi_k'(\theta) = \lim_{\theta \to \pi} \Psi_k(\theta) =
\lim_{\theta \to \pi} \sin \theta \Psi_k'(\theta) = 0.  \label{X}
\end{equation}

{\bf Case $k = 0$:} It follows from the divergence-free condition
(\ref{2.9}) that
$$
U_0 = \frac{\alpha}{\sin \theta},
$$
where $\alpha \in \mathbb{R}$. This solution resembles the neutral
eigenmode generated by the arbitrary constant $\alpha$ in the
stationary solution (\ref{exact-solution-2D}). Since the eigenmode
violates the boundary conditions (\ref{X0}) on $S_0$ and has pole
singularities on $S$, we set $\alpha = 0$. In this case, the first
equation (\ref{2.7}) admits a solution $Q_0 = \beta$, where $\beta
\in \mathbb{R}$. It is also a neutral eigenmode generated by the
arbitrary constant $\beta$ in the stationary solution
(\ref{exact-solution-2D}). Since it is a trivial eigenmode (the
pressure term is defined with accuracy to an addition of an
arbitrary constant), we can set $\beta = 0$.

When $\alpha = \beta = 0$, the representation for $U_0 = Q_0 = 0$
matches the previous representation for $U_k$ and $Q_k$ with $k =
0$. Using the representation (\ref{stream-function-representation}),
we introduce $V_0 = -\Psi_0'(\theta)$ and rewrite the second
equation (\ref{2.8}) as follows:
\begin{equation}
\label{zero-equation} \frac{d}{d\theta}\left( \nu \Delta_0 \Psi_0 -
\lambda \Psi_0 \right) - \frac{1}{\sin \theta} \Delta_0 \Psi_0 = 0,
\end{equation}
where $\Delta_0$ is defined by (\ref{Laplacian-k}) with $k = 0$.
Letting $\Phi_0 = \Delta_0 \Psi_0$ and taking one more derivative
in $\theta$, one can convert the non-trivial equation
(\ref{zero-equation}) to the previous form
(\ref{eigenvalue-equation}) with $k = 0$. Therefore, all solutions
of (\ref{zero-equation}) are also solutions of
(\ref{eigenvalue-equation}) with $k = 0$, while the converse
statement is not true. It follows from (\ref{bc-dirichlet-0}) and
(\ref{bc-dirichlet}) that the stream function $\Psi_0(\theta)$
satisfies the Neumann boundary conditions
\begin{equation}
\Psi_0'(\theta_0) = \Psi_0'(\pi-\theta_0) = 0 \label{X0-zero}
\end{equation}
in the case of $S_0$ and the boundary conditions
\begin{equation}
\lim_{\theta \to 0} \sin \theta \Psi_0'(\theta) = \lim_{\theta \to
\pi} \sin \theta \Psi_0'(\theta) = 0  \label{X-zero}
\end{equation}
in the case of $S$. Stability analysis of the linearized system
(\ref{eigenvalue-equation})--(\ref{second-equation}) with $k \neq
0$ is developed separately from that of the linearized equation
(\ref{zero-equation}) with $k = 0$. Our main results on
eigenvalues of the linearized systems
(\ref{eigenvalue-equation})--(\ref{second-equation}) and
(\ref{zero-equation}) are summarized in Table 1. The remainder of
this article is devoted to the proofs and numerical verifications
of results described in Table 1.

\vspace{0.1cm}

\begin{tabular}{|p{1.5cm}|p{2cm}|p{2cm}|p{4cm}|p{4cm}|}
\hline
Index $k$ & Viscosity $\nu$ & Cut-off $\theta_0$ & eigenvalues & results \\
\hline
$k \neq 0$ & $0 < \nu \leq \infty$ & $\theta_0 = 0$ & real negative & Proposition 2 \\
\hline $k \neq 0$ & $\nu = \infty$ & $0 < \theta_0 <
\frac{\pi}{2}$ & real
negative & Propositions 3 and 4 \\
\hline $k \neq 0$ & $0 < \nu < \infty$ & $0 < \theta_0 <
\frac{\pi}{2}$ & real or complex & Section 5 \\
\hline $k = 0$ & $0 < \nu \leq \infty$ & $\theta_0 = 0$ & real
negative or absent & Proposition 8 \\
\hline $k = 0$ & $0 < \nu \leq \infty$ & $0 < \theta_0 < \frac{\pi}{2}$ & real negative & Propositions 9 and 10 \\
\hline $k = 0$ & $0 < \nu < \infty$ & $0 < \theta_0 < \frac{\pi}{2}$ & real negative & Section 6 \\
\hline
\end{tabular}

{\bf Table 1:} Summary of main results.

\section{Stability analysis for $k \neq 0$}

We rewrite the coupled system
(\ref{eigenvalue-equation})--(\ref{second-equation}) for
$(\Psi_k,\Phi_k)$ by using the variable $x = \cos \theta$:
\begin{equation}
\label{stability-problem} L_k \Psi_k = \Phi_k, \qquad L_k \Phi_k +
\epsilon \Phi_k' = \mu \Phi_k,
\end{equation}
where $\epsilon = 1/\nu$ is the Reynolds number of the basic flow,
$\mu = \lambda/\nu$ is a rescaled eigenvalue, and $L_k$ is the
Sturm--Liouville operator for associated Legendre functions
\begin{equation}
L_k = \frac{d}{d x} \left[ (1-x^2) \frac{d}{dx} \right] -
\frac{k^2}{1-x^2}. \label{Sturm-Liouville}
\end{equation}
The system (\ref{stability-problem}) is defined on the symmetric
interval $-x_0 \leq x \leq x_0$, where $x_0 = \cos \theta_0$. The
spherical layer $S_0$ corresponds to the case $0 < x_0 < 1$, while
the complete sphere $S$ corresponds to the limit $x_0 \to 1$. In
the latter case, the interval $x \in [-1,1]$ connects two singular
points $x = \pm 1$ of the Sturm--Liouville operator
(\ref{Sturm-Liouville}). The case $\epsilon = 0$ corresponds to
the infinitely viscous fluid, while the case $\epsilon = \infty$
corresponds to the inviscous fluid.

Using the representation (\ref{stream-function-representation})
and the transformation $x = \cos \theta$ with $\Psi_k'(\theta) = -
\sqrt{1 - x^2} \Psi_k'(x)$, we rewrite the condition
(\ref{square-integrable}) as the norm on function space ${\cal
H}_k$, which is used throughout our work:
\begin{equation}
\label{energy-norm} \| \Psi_k \|^2_{{\cal H}_k} = \int_{-x_0}^{x_0}
\left[ (1-x^2) |\Psi_k'(x)|^2 + \frac{k^2}{1 - x^2} |\Psi_k(x)|^2
\right] dx < \infty.
\end{equation}
We shall denote ${\cal H}_k([-x_0,x_0])$ when $0 < x_0 < 1$ and
${\cal H}_k([-1,1])$ when $x_0 = 1$. When $0 < x_0 < 1$, the
linearized system (\ref{stability-problem}) is defined on function
space
\begin{equation}
X_0 = \left\{ \Psi_k \in {\cal H}_k([-x_0,x_0]) : \quad \Psi_k(\pm
x_0) = \Psi_k'(\pm x_0) = 0 \right\}, \label{bc1}
\end{equation}
where the boundary conditions (\ref{X0}) are taken into account.
When $x_0 = 1$, the linearized system (\ref{stability-problem}) is
defined in function space
\begin{equation}
\label{bc2} X = \left\{ \Psi_k \in {\cal H}_k([-1,1]) : \quad
\lim_{x \to \pm 1} \Psi_k(x) = \lim_{x \to \pm 1} (1-x^2) \Psi'_k(x)
= 0 \right\},
\end{equation}
where the boundary conditions (\ref{X}) are taken into account. We
note that the boundary conditions in the definition of $X$ are
redundant, since the norm (\ref{energy-norm}) is finite on $x \in
[-1,1]$ only if the boundary conditions in (\ref{bc2}) are
satisfied. Nevertheless, we write these redundant boundary
conditions according to the standard formalism of the singular
Sturm--Lioville problems \cite{Strauss}.

The Sturm--Liouville operator $L_k$ in (\ref{Sturm-Liouville}) is
self-adjoint with respect to the boundary conditions in $X_0$ and
$X$, such that $(\Psi_k, L_k \Psi_k) = -\| \Psi_k \|^2_{{\cal
H}_k} < 0$ is finite and real-valued for $\Psi_k \in {\cal
H}_k([-x_0,x_0])$. Therefore, the kernel of $L_k$ is empty in
$X_0$ and $X$. Because the smallest eigenvalue of $L_k$ is bounded
away zero, the operator $L_k$ is invertible and ${\rm range}(L_k)$
is dense in the space of square integrable functions on $x \in
[-x_0,x_0]$ for any $0 < x_0 \leq 1$. Therefore, as it follows
from the first equation of the system (\ref{stability-problem}),
the component $\Phi \in {\rm range}(L_k)$ is square integrable on
$x \in [-x_0,x_0]$ but it does not satisfy any specific boundary
conditions at the end points $x = \pm x_0$.

The eigenvalue problem (\ref{stability-problem}) in $X_0$ and $X$
has two continuous parameters $0 < x_0 \leq 1$ and $\epsilon \geq 0$
and one integer parameter $k \in \mathbb{Z} \backslash \{0\}$, while
$(\mu,\Psi_k)$ is the eigenvalue-eigenfunction pair that defines
spectral stability of the stationary flow. The following results
characterize the spectrum of the eigenvalue problem in the cases:
(i) $x_0 = 1$ and $\epsilon \geq 0$; (ii) $0 < x_0 < 1$ and
$\epsilon = 0$; and (iii) in the limit $x_0 \to 1$ when $\epsilon =
0$. Based on these results, we prove the following theorem:

\begin{theorem}
When $x_0 = 1$ and $\epsilon \geq 0$ or $0 < x_0 \leq 1$ and
$\epsilon = 0$, the stationary flow (\ref{exact-solution-2D}) is
asymptotically stable with respect to symmetry-breaking
perturbations in the sense that the spectrum of the linearized
problem (\ref{stability-problem}) in $X_0$ or $X$ consists of a
set of isolated eigenvalues $\mu$ of finite multiplicities, where
$\mu \in \mathbb{R}_-$ is bounded away from zero.
\label{theorem-1}
\end{theorem}

The proof of theorem consists of the proofs of three individual
propositions.

\begin{proposition}
\label{proposition-x0-1} A complete spectrum of the eigenvalue
problem (\ref{stability-problem}) with $x_0 = 1$ and $\epsilon \geq
0$ in $X$ consists of simple isolated eigenvalues at $\mu = \mu_n$,
\begin{equation}
\label{eigenvalues-1} \mu_n = -s_n(s_n+1), \qquad s_n = \sigma + n,
\end{equation}
where $\sigma = \sqrt{k^2 + \epsilon^2/4} > 0$ and $n \geq 0$ is
integer.
\end{proposition}

\begin{proof}
Let $\mu = -s(s+1)$ and
\begin{equation}
\label{transformation-1} \Phi_k(x) = \left( \frac{1-x}{1+x}
\right)^{\epsilon/4} \varphi(x).
\end{equation}
The second equation of the system (\ref{stability-problem})
transforms to the associated Legendre equation
\begin{equation}
\label{associated-Legendre} \frac{d}{d x} \left[ (1-x^2) \frac{d
\varphi}{dx} \right] - \frac{\sigma^2}{1-x^2} \varphi + s(s+1)
\varphi = 0, \qquad -1 < x < 1,
\end{equation}
where $\sigma = \sqrt{k^2 + \epsilon^2/4} > 0$. Since the linear ODE
(\ref{associated-Legendre}) has no singular points on $-1 < x < 1$,
there exists a set of two linearly independent twice continuously
differentiable solutions in any compact subset of $x \in (-1,1)$
\cite{CL55}. Singularity analysis of the ODE
(\ref{associated-Legendre}) as $x \to \pm 1$ shows that the solution
$\varphi(x)$ either have a singular (unbounded) behavior like $(1
\mp x)^{-\sigma/2}$ as $x \to \pm 1$ or a regular (vanishing)
behavior like $(1 \mp x)^{\sigma/2}$ as $x \to \pm 1$.

Let $\varphi(x)$ be a regular solution of
(\ref{associated-Legendre}) on $x \in [-1,1]$, such that $\varphi(x)
\sim (1\mp x)^{\sigma/2}$ and $\Phi_k(x) \sim (1 \mp x)^{\pm
\epsilon/4 + \sigma/2}$ as $x \to \pm 1$. Since the Sturm--Liouville
operator $L_k$ is invertible on $\Phi_k \in L^2([-1,1])$ for $k \neq
0$, the first equation of the system (\ref{stability-problem})
admits a solution $\Psi_k(x)$ that behaves like $(1 \mp x)^{1 \pm
\epsilon/4 + \sigma/2}$ as $x \to \pm 1$. Since $\pm \epsilon +
\sqrt{\epsilon^2 + 4 k^2} \geq 0$ for any $k \in \mathbb{Z}$ and
$\epsilon \geq 0$, the function $\Phi_k(x)$ is bounded on $x \in
[-1,1]$ such that $\Phi_k \in L^2([-1,1])$ while the function
$\Psi_k(x)$ belongs to the function space $X$ in (\ref{bc2}).
Therefore, if $\varphi(x)$ is a regular solution of
(\ref{associated-Legendre}), then $\Psi_k(x)$ is an eigenfunction of
the eigenvalue problem (\ref{stability-problem}) in $X$.

Let $\varphi(x)$ be a singular solution of
(\ref{associated-Legendre}), such that $\varphi(x) \sim (1 \mp
x)^{-\sigma/2}$ and $\Phi_k(x) \sim (1 \mp x)^{\pm \epsilon/4 -
\sigma/2}$ in at least one limit $x \to \pm 1$. Since $\pm
\epsilon - \sqrt{\epsilon^2 + 4 k^2} \leq - 2 |k| \leq -2$ for $k
\neq 0$ and $\epsilon \geq 0$, the function $\Phi_k$ does not
belong to $L^2([-1,1])$ and hence $\Psi_k(x)$ can not be in $X$.
By Theorem 10 on p.1441 in \cite{DS}, the essential spectrum of
the formally self-adjoint operator (\ref{associated-Legendre}) is
void. Therefore, the complete spectrum of the linearized system
(\ref{stability-problem}) in $X$ consists of isolated eigenvalues
$\mu$, which correspond to {\em regular} solutions $\varphi(x)$ of
the associated Legendre equation (\ref{associated-Legendre}).

Let $\varphi(x)$ be a regular solution of
(\ref{associated-Legendre}) and write $\varphi(x) =
(1-x^2)^{\sigma/2} F(x)$, where $F(x)$ is bounded as $x \to \pm 1$.
This substitution transforms the associated Legendre equation
(\ref{associated-Legendre}) to the hypergeometric equation
\begin{equation}
\label{hypergeometric-equation} z(1-z) F''(z) + \left( \gamma -
(\alpha + \beta + 1)z\right) F'(z) - \alpha \beta F(z) = 0,
\end{equation}
where
\begin{equation}
\label{hypergeometric-parametrization} z = \frac{1-x}{2}, \quad
\alpha = \sigma - s, \quad \beta = \sigma + s + 1, \quad \gamma =
\sigma + 1.
\end{equation}
The only solution of the ODE (\ref{hypergeometric-equation}) which
is bounded as $x \to 1$ ($z \to 0$) is the hypergeometric function
$F(z;\alpha,\beta,\gamma)$, which admits the power series at $z = 0$
(see 9.100 on p. 995 in \cite{GR}):
\begin{equation}
\label{hypergeometric} F(z;\alpha,\beta,\gamma) = 1 + \frac{\alpha
\beta}{\gamma 1!} z + \frac{\alpha (\alpha+1) \beta (\beta +
1)}{\gamma (\gamma + 1) 2!} z^2 + ...
\end{equation}
The hypergeometric series (\ref{hypergeometric}) converges for $|z|
< 1$ but it diverges as $z \to 1$ since $\alpha + \beta - \gamma =
\sigma > 0$ unless the truncation of the power series to a
polynomial in $z$ occurs (see 9.101--9.102 on p. 995 in \cite{GR}).
The latter case is the only case when the solution of the ODE
(\ref{hypergeometric-equation}) is bounded in both limits $x \to 1$
($z \to 0$) and $x \to -1$ ($z \to 1$). It is easy to see that the
truncation occurs when either $\alpha = -n$ or $\beta = -m$ with
non-negative integers $n$ and $m$. The two cases are in fact
equivalent to each other since $\mu = -s(s+1) = (\alpha - \sigma)
(\beta - \sigma)$ and $\alpha + \beta = 1 + 2 \sigma$. Let $\alpha =
-n$, such that $s = \sigma + n$, $\beta = 2\sigma + n + 1$ and
$\gamma = \sigma + 1$. In this case, the function
$F\left(z;-n,n+1+2\sigma,1+\sigma\right) \equiv F_n(x)$ is a
polynomial of degree $n$, e.g.
\begin{equation}
F_0 = 1, \quad F_1 = x, \quad F_2 = \frac{(2\sigma+3) x^2 -
1}{2(1+\sigma)}, \quad F_3 = \frac{(5 + 2 \sigma) x^3 - 3
x}{2(1+\sigma)}, \label{polynomials-Legendre}
\end{equation}
while the simple eigenvalues $\mu = \mu_n$ are given by the
expression (\ref{eigenvalues-1}). When $\sigma = 0$, polynomials
$F_n(x)$ coincide with the Legendre polynomials $P_n(x)$ in 8.91 on
p. 973 of \cite{GR}.
\end{proof}

\begin{proposition}
\label{proposition-epsilon-0} A complete spectrum of the eigenvalue
problem (\ref{stability-problem}) with $0 < x_0 < 1$ and $\epsilon =
0$ in $X_0$ consists of isolated eigenvalues $\mu$, which are (i)
real and strictly negative and (ii) either simple or double with
linearly independent eigenfunctions.
\end{proposition}

\begin{proof}
We first show that no zero eigenvalue $\mu = 0$ exists in the
eigenvalue problem (\ref{stability-problem}) with $0 < x_0 < 1$ and
$\epsilon = 0$ in $X_0$. Let $\Psi_k(x)$ be a $C^4([-x_0,x_0])$
solution of the fourth-order ODE $L_k^2 \Psi_k = 0$ in function
space $X_0$. Then,
$$
(\Psi_k,L_k^2 \Psi_k) = (1-x^2) \left[ \Psi_k (L_k \Psi_k)' -
\Psi_k' (L_k \Psi_k) \right] |_{x=-x_0}^{x=x_0} + (L_k \Psi_k, L_k
\Psi_k) = (L_k \Psi_k, L_k \Psi_k).
$$
Therefore, $\Psi_k(x)$ is in fact the solution of the second-order
ODE $L_k \Psi_k = 0$. The boundary conditions in $X_0$ admit the
only solution $\Psi_k(x) \equiv 0$, such that the eigenvalue problem
(\ref{stability-problem}) contains no eigenvalue $\mu = 0$ in $X_0$.

When $\mu \neq 0$ and $\epsilon = 0$, the system
(\ref{stability-problem}) admits a general solution in the form
$$
\Psi_k(x) = \frac{\phi(x)}{\mu} + \psi(x), \qquad \Phi_k(x) =
\phi(x),
$$
where $\psi(x)$ and $\phi(x)$ are general solutions of the
homogeneous second-order ODEs
$$
L_k \psi = 0, \qquad L_k \phi = \mu \phi.
$$
Since the operator $L_k$ is invariant with respect to the inversion
symmetry $x \mapsto - x$, each homogeneous second-order ODE has
linearly independent symmetric (even) and anti-symmetric (odd)
solutions denoted by subscripts $+$ and $-$ respectively. Therefore,
we obtain the decomposition
\begin{eqnarray*}
\Psi_k(x) & = & d_+ \frac{\phi_+(x)}{\mu} + c_+ \psi_+(x) + d_-
\frac{\phi_-(x)}{\mu} + c_- \psi_-(x), \\
\Phi_k(x) & = & d_+ \phi_+(x) + d_- \phi_-(x),
\end{eqnarray*}
where $(c_+,c_-,d_+,d_-)$ are constants and the functions
$\phi_{\pm}(x)$ and $\psi_{\pm}(x)$ are uniquely normalized by the
initial values at $x = 0$ (e.g. $\phi_+(0) = 1$, $\phi_+'(0) = 0$
and $\phi_-(0) = 0$, $\phi_-'(0) = 1$). We note that either
$\psi_{\pm}(x_0) \neq 0$ or $\psi_{\pm}'(x_0) \neq 0$ (since
$\psi_{\pm}(x) \equiv 0$ otherwise). By using the boundary
conditions in (\ref{bc1}), we decompose the boundary-value problems
into two uncoupled systems with
$$
d_{\pm} \phi_{\pm}(x_0) + \mu c_{\pm} \psi_{\pm}(x_0) = 0, \qquad
d_{\pm} \phi_{\pm}'(x_0) + \mu c_{\pm} \psi_{\pm}'(x_0) = 0,
$$
such that a non-zero solution for $(c_+,c_-,d_+,d_-)$ exists
provided
$$
\phi_{\pm}'(x_0) \psi_{\pm}(x_0) = \phi_{\pm}(x_0) \psi'_{\pm}(x_0).
$$
The functions $\psi_{\pm}(x)$ are independent of $\mu$, while
$\phi_{\pm}(x)$ depend on $\mu$. We have thus obtained that the
functions $\phi_{\pm}(x)$ solve the {\em closed} eigenvalue problem
\begin{equation}
\label{equivalent-problem-1} L_k \phi_{\pm} = \mu \phi_{\pm}, \qquad
-x_0 \leq x \leq x_0,
\end{equation}
defined on the function space
\begin{equation}
\label{equivalent-problem-2} H_0 = \left\{ \phi_{\pm} \in {\cal
H}_k([-x_0,x_0]) : \;\; \psi_{\pm}(x_0) \phi'_{\pm}(x_0) -
\psi_{\pm}'(x_0) \phi_{\pm}(x_0) = 0, \;\; \phi_{\pm}(-x) = \pm
\phi_{\pm}(x) \right\}.
\end{equation}
The $\mu$-independent boundary values in
(\ref{equivalent-problem-2}) are Robin boundary conditions when
$\psi_{\pm}(x_0)$ and $\psi_{\pm}'(x_0)$ are both non-zero,
Dirichlet boundary conditions when $\psi_{\pm}(x_0) = 0$ and Neumann
boundary conditions when $\psi'_{\pm}(x_0) = 0$. The associated
Legendre operator $L_k$ is self-adjoint in $H_0$ with respect to any
of these boundary conditions \cite{Strauss}. Therefore, all
eigenvalues $\mu$ of the eigenvalue problem
(\ref{equivalent-problem-1}) in $H_0$ are real-valued and isolated,
while the corresponding eigenfunctions $\phi_{\pm}(x)$ are
real-valued. Moreover, all eigenvalues of
(\ref{equivalent-problem-1}) are simple since the Wronskian of any
two solutions of (\ref{equivalent-problem-1}) with boundary
conditions in (\ref{equivalent-problem-2}) is zero. Since $\Psi_k
\in X_0$ and $\Phi_k \in H_0$, we obtain that
$$
(\phi,\phi) = (L_k \Psi_k, \phi) = (\Psi_k, L_k \phi) = \mu
(\Psi_k,\phi) = (\phi,\phi) + \mu (\psi,\phi),
$$
such that $(\psi,\phi) = 0$ for $\mu \neq 0$. By using the above
identity, we obtain that
\begin{equation}
\label{Green-identity} \frac{1}{\mu} (\phi,\phi) = (\Psi_k,\phi) =
(\Psi_k, L_k \Psi_k) = - \| \Psi_k \|^2_{{\cal H}_k} < 0,
\end{equation}
such that $\mu < 0$ for {\em each} eigenvalue with $\Psi_k \neq 0$
and $\phi \neq 0$. By construction, eigenvalues are at most
double. The case of double eigenvalues corresponds to the
situation when the eigenvalue problems
(\ref{equivalent-problem-1})--(\ref{equivalent-problem-2}) admit
two linearly independent (even and odd) eigenfunctions for the
same value of $\mu$.
\end{proof}

\begin{proposition}
Let $\{ \mu_n \}_{n \geq 0}$ be isolated eigenvalues of the
eigenvalue problem (\ref{stability-problem}) in $X_0$ with $0 < x_0
< 1$ and $\epsilon = 0$ ordered as
$$
\mu_0 \geq \mu_1 \geq ... \geq \mu_n \geq ...
$$
Then,
$$
\lim_{x_0 \to 1} \mu_n = - s_n(s_n+1), \qquad s_n = |k| + n,
$$
where $n \geq 0$. \label{proposition-continuity}
\end{proposition}

\begin{proof}
Consider even and odd solutions of the second-order ODE $L_k
\psi_{\pm} = 0$ in the limit $x_0 \to 1$. Since the kernel of $L_k$
admits no eigenfunctions in ${\cal H}_k$ for $k \neq 0$ and $0 < x_0
\leq 1$, the solutions $\psi_{\pm}(x_0)$ must diverge as $x_0 \to
1$. Singularity analysis as $x \to \pm 1$ suggests that the solution
$\psi_{\pm}(x)$ grows like $(1 \mp x)^{-|k|/2}$ as $x \to \pm 1$,
such that $\lim\limits_{x_0 \to 1} \psi_{\pm}(x_0)/\psi'_{\pm}(x_0)
= 0$. Therefore, eigenfunctions $\phi_{\pm}(x)$ of the auxiliary
eigenvalue problem (\ref{equivalent-problem-1}) for $0 < x_0 < 1$
satisfy in the limit $x_0 \to 1$ the singular eigenvalue problem
\begin{equation}
\label{equivalent-problem-3} L_k \phi_{\pm} = \mu \phi_{\pm}, \qquad
-1 < x < 1
\end{equation}
defined on the function space
\begin{equation}
\label{equivalent-problem-4} H = \left\{ \phi_{\pm} \in {\cal
H}_k([-1,1]) : \quad \lim_{x \to \pm 1} \phi_{\pm}(x) = \lim_{x \to
\pm 1} (1-x^2) \phi_{\pm}'(x) = 0 \right\}.
\end{equation}
Again, the boundary conditions in $H$ are redundant due to
convergence of the integral in ${\cal H}_k([-1,1])$. A complete
spectrum of the eigenvalue problem
(\ref{equivalent-problem-3})--(\ref{equivalent-problem-4}) is
constructed in the proof of Proposition \ref{proposition-x0-1}:
eigenvalues are given by (\ref{eigenvalues-1}) with $\epsilon = 0$
and eigenfunctions are $\phi_{\pm}(x) = (1-x^2)^{|k|/2} F_n(x)$,
where $F_n(x)$ are associated Legendre polynomials
(\ref{polynomials-Legendre}) with $\sigma = |k|$. Convergence and
uniqueness of continuations from eigenvalues of
(\ref{equivalent-problem-1}) in $H_0$ for $x_0 < 1$ to eigenvalues
of (\ref{equivalent-problem-3}) in $H$ for $x_0 = 1$ is proved in
two steps. Theorem 5.3 of \cite{zettle} guarantees convergence and
uniqueness of continuations from the singular Sturm--Liouville
problem (\ref{equivalent-problem-3}) in $H$ to the regular
Dirichlet problem for the Sturm--Liouville operator
(\ref{equivalent-problem-1}) on $-x_0 \leq x \leq x_0$. The
Dirichlet problem is generally different from the Robin
boundary-value problem in $H_0$ by the terms $\psi_{\pm}(x_0)
\phi_{\pm}'(\pm x_0)/\psi_{\pm}'(x_0)$ in the boundary conditions
in $H_0$. However, these terms are small in the limit $x_0 \to 1$.
Unique continuation of simple eigenvalues of the Dirichlet problem
to the simple eigenvalues of the Robin problem (separately for
$\phi_+(x)$ and $\phi_-(x)$) follows by standard perturbation
theory of eigenvalues of self-adjoint Sturm--Liouville operators
in Lemma VIII 1.24 of \cite{Kato}.
\end{proof}

\begin{remark}
{\rm Theorem \ref{theorem-1} does not cover the case $0 < x_0 < 1$
and $\epsilon > 0$. Eigenvalues of the linearized problem
(\ref{stability-problem}) in this case will be computed in Section 5
numerically. } \label{remark-k}
\end{remark}

\section{Stability analysis for $k = 0$}

We rewrite the linearized equation (\ref{zero-equation}) in the
variable $x = \cos \theta$:
\begin{equation}
\label{zero-equation-stability} L_0 \Psi_0 = \Phi_0, \qquad \Phi'_0
+ \frac{\epsilon}{1 - x^2} \Phi_0 = \mu \Psi_0'.
\end{equation}
where $L_0$ is the Sturm--Liouville operator for Legendre functions
\begin{equation}
L_0 = \frac{d}{d x} \left[ (1-x^2) \frac{d}{dx} \right]
\label{Sturm-Liouville-0}
\end{equation}
and $\Phi_0(x)$ is introduced similarly to the system
(\ref{stability-problem}). Incorporating the boundary conditions
(\ref{X0-zero}) and (\ref{X-zero}) in new variables, we introduce
the function spaces $X_0$ and $X$ for the eigenvalue problem
(\ref{zero-equation-stability}). When $0 < x_0 < 1$, the function
space $X_0$ is
\begin{equation}
X_0 = \left\{ \Psi_0 \in {\cal H}_0([-x_0,x_0]) : \quad \Psi_0'(\pm
x_0) = 0 \right\}. \label{bc1-zero}
\end{equation}
When $x_0 = 1$, the function space $X$ is
\begin{equation}
X = \left\{ \Psi_0 \in {\cal H}_0([-1,1]) : \quad \lim_{x \to \pm 1}
(1-x^2) \Psi_0'(x) = 0 \right\}, \label{bc2-zero}
\end{equation}
where the boundary conditions are redundant due to convergence of
the integral in ${\cal H}_0([-1,1])$. No boundary conditions on
$\Psi_0(x)$ are set at $x = \pm x_0$. Moreover, the system
(\ref{zero-equation-stability}) defines the function $\Psi_0(x)$
up to an arbitrary additive constant. Therefore, the constant
function $\Psi_0(x) \equiv {\rm const}$ is always an eigenfunction
of the system (\ref{zero-equation-stability}) with $\Phi_0(x)
\equiv 0$.

\begin{lemma}
The eigenvalue $\mu = 0$ of the linearized system
(\ref{zero-equation-stability}) in either $X_0$ or $X$ is
algebraically and geometrically simple. \label{lemma-zero}
\end{lemma}

\begin{proof}
Integrating the first equation in the system
(\ref{zero-equation-stability}) on $x \in [-x_0,x_0]$ for
$\Psi_0(x)$ in either $X_0$ or $X$, we obtain the Fredholm
Alternative condition
\begin{equation}
\label{Fredholm} \int_{-x_0}^{x_0} \Phi_0(x) dx = 0,
\end{equation}
where $0 < x_0 \leq 1$. Integrating the second equation in the
system (\ref{zero-equation-stability}), we obtain a general solution
for $\mu = 0$:
$$
\Phi_0 = c_0 \left(\frac{1 - x}{1+x}\right)^{\epsilon/2},
$$
where $c_0$ is constant. Since $\Phi_0(x)$ does not satisfy the
Fredholm Alternative condition (\ref{Fredholm}), we have to set $c_0
= 0$. Then, $\Psi_0(x)$ satisfies the second-order ODE $L_0 \Psi_0 =
0$, which admits only one eigenfunction $\Psi_0(x) \equiv {\rm
const}$ in either $X_0$ or $X$. Similarly one can prove that the
Jordan block of the zero eigenvalue with the eigenfunction
$\Psi_0(x) \equiv {\rm const}$ and $\Phi_0(x) \equiv 0$ is of the
length one.
\end{proof}

We will extend results of Section 3 to the linearized problem
(\ref{zero-equation-stability}) with $\mu \neq 0$ in $X_0$ and
$X$. Neglecting the only zero eigenvalue $\mu = 0$ with the
trivial eigenfunction $\Psi_0(x) \equiv {\rm const}$, we prove the
following theorem.

\begin{theorem}
The stationary flow (\ref{exact-solution-2D}) is asymptotically
stable with respect to symmetry-preserving perturbations in the
sense that all eigenvalues $\mu$ (excluding the trivial zero) of the
linearized problem (\ref{zero-equation-stability}) with $0 < x_0
\leq 1$ and $\epsilon \geq 0$ in $X_0$ or $X$ are real and strictly
negative.   \label{theorem-2}
\end{theorem}

In order to develop analysis of eigenvalues for $\mu \neq 0$, we
shall use two equivalent reformulations of the third-order ODE
system (\ref{zero-equation-stability}) as the second-order
eigenvalue problems associated with formally self-adjoint operators.
In the first reformulation, we exclude $\mu \Psi_0'(x)$ from the
system (\ref{zero-equation-stability}) and find a closed equation
for $\Phi_0(x)$,
\begin{equation}
\label{stability-zero-problem} L_0 \Phi_0 + \epsilon \Phi_0' = \mu
\Phi_0.
\end{equation}
By introducing new dependent variable $\varphi(x)$ via
\begin{equation}
\label{transformation-2} \Phi_0(x) = \left( \frac{1-x}{1+x}
\right)^{\epsilon/4} \varphi(x),
\end{equation}
the linearized equation (\ref{stability-zero-problem}) is
transformed to the self-adjoint form given by the associated
Legendre equation
\begin{equation}
\label{associated-Legendre-zero} \frac{d}{d x} \left[ (1-x^2)
\frac{d \varphi}{dx} \right] - \frac{\epsilon^2}{4(1-x^2)} \varphi =
\mu \varphi, \qquad -x_0 < x < x_0.
\end{equation}
By using the second equation of the system
(\ref{zero-equation-stability}), we obtain the first-order ODE for
the function $\Psi_0(x)$:
\begin{equation}
\label{Psi-0-coupling} \mu \Psi_0'(x) = \left( \frac{1 -
x}{1+x}\right)^{\epsilon/4} \left( \frac{d \varphi}{dx} +
\frac{\epsilon}{2(1-x^2)} \varphi \right).
\end{equation}
While the linearized equation (\ref{stability-zero-problem})
coincides with the second equation of the system
(\ref{stability-problem}) for $k = 0$, the present role of this
equation is different. In order to find $\Psi_0(x)$ from a
solution $\Phi_0(x)$ of the closed equation
(\ref{stability-zero-problem}), we can solve the first-order ODE
(\ref{Psi-0-coupling}) in either $X_0$ or $X$ with $\mu \neq 0$.
Therefore, as opposed to the case $k \neq 0$, we do not have to
solve the first equation of the system
(\ref{zero-equation-stability}) and the Fredholm Alternative
condition (\ref{Fredholm}) can be ignored in this approach.

In the second reformulation of the third-order ODE system
(\ref{zero-equation-stability}), we introduce a new dependent
variable $\chi(x)$ via
\begin{equation}
\label{transformation-3} \Psi_0'(x) = \left( \frac{1-x}{1+x}
\right)^{\epsilon/4} \frac{\chi(x)}{\sqrt{1-x^2}}.
\end{equation}
By using the first equation of the system
(\ref{zero-equation-stability}), we express the function $\Phi_0(x)$
in terms of $\chi(x)$:
\begin{equation}
\label{Phi-0-coupling} \Phi_0(x) = \left( \frac{1 -
x}{1+x}\right)^{\epsilon/4} \left( \sqrt{1-x^2} \frac{d \chi}{dx} -
\frac{\epsilon + 2x}{2\sqrt{1-x^2}} \chi \right).
\end{equation}
The second equation of the system (\ref{zero-equation-stability})
transforms then to the self-adjoint form:
\begin{equation}
\label{self-adjoint-form} \frac{d}{d x} \left[ (1-x^2) \frac{d
\chi}{dx} \right] - \frac{\epsilon^2 + 4 + 4 \epsilon x}{4(1-x^2)}
\chi = \mu \chi, \qquad -x_0 < x < x_0.
\end{equation}
Although the second-order ODE (\ref{self-adjoint-form}) is more
complicated than the associated Legendre equation
(\ref{associated-Legendre-zero}), the eigenfunction $\chi(x)$ is
related to the function $\Psi_0(x)$ better than the eigenfunction
$\varphi(x)$. In particular, when $x_0 = 1$ and $\Psi_0 \in X$, the
eigenfunction $\chi(x)$ satisfies the conditions:
\begin{equation}
\label{conditions-chi-X} \int_{-1}^1 \left( \frac{1-x}{1+x}
\right)^{\epsilon/2} \chi^2(x) dx < \infty, \qquad \lim_{x \to \pm
1} \left( \frac{1-x}{1+x} \right)^{\epsilon/4} \sqrt{1-x^2} \chi(x)
= 0.
\end{equation}
When $0 < x_0 < 1$ and $\Psi_0 \in X_0$, the eigenfunction $\chi(x)$
is any classical solution of the second-order ODE
(\ref{self-adjoint-form}) on $x \in [-x_0,x_0]$ with the Dirichlet
boundary conditions $\chi(\pm x_0) = 0$. There exists a pair of
Darboux-Backlund transformations between the Sturm--Liouville
problems (\ref{associated-Legendre-zero}) and
(\ref{self-adjoint-form}):
\begin{eqnarray}
\label{backlund-1} \varphi(x) & = & \sqrt{1 - x^2} \chi'(x) -
\frac{\epsilon + 2 x}{2 \sqrt{1 - x^2}} \chi(x), \\
\label{backlund-2} \mu \chi(x) & = & \sqrt{1 - x^2} \varphi'(x) +
\frac{\epsilon}{2 \sqrt{1 - x^2}} \varphi(x),
\end{eqnarray}
where $\mu \neq 0$ is assumed. By the Friedrichs' theorems (see,
e.g. Theorem 10 on p.1441 or Theorem 67 on p. 1501 of \cite{DS}),
the essential spectrum of the formally self-adjoint operators
(\ref{associated-Legendre-zero}) and (\ref{self-adjoint-form}) is
void. Therefore, the spectrum of these operators consists of a
sequence of isolated eigenvalues of finite multiplicities, which
we identify in three individual propositions.

\begin{proposition}
\label{proposition-x0-1-k-0} A complete spectrum of the eigenvalue
problem (\ref{zero-equation-stability}) with $x_0 = 1$ and $0 \leq
\epsilon < 2$ in $X$ consists of simple isolated eigenvalues at $\mu
= \mu_n$, where
\begin{equation}
\label{eigenvalues-1-zero} \mu_n = - n(n+1), \qquad n \geq 0.
\end{equation}
No non-zero eigenvalues of the eigenvalue problem
(\ref{zero-equation-stability}) with $x_0 = 1$ and $\epsilon \geq 2$
exists in $X$.
\end{proposition}

\begin{proof}
Let $\mu = -s(s+1) \neq 0$ and $\varphi(x) = (1-x^2)^{\epsilon/4}
F(x)$ and consider the associated Legendre equation
(\ref{associated-Legendre-zero}) with $x_0 = 1$. Then, the
function $F(x)$ satisfies the hypergeometric equation
(\ref{hypergeometric-equation}) under parametrization
(\ref{hypergeometric-parametrization}) with $\sigma = \epsilon/2$.
In order to identify solutions $F(x)$ of the hypergeometric
equations in the function space $\Psi_0 \in X$, we shall rewrite
the relation (\ref{Psi-0-coupling}) as follows:
\begin{equation}
\label{Psi-0-equation} -s(s+1) \Psi_0'(x) = (1 - x)^{\epsilon/2}
\left( F'(x) + \frac{\epsilon}{2(1+x)} F(x) \right).
\end{equation}
Also recall that $\Phi_0(x) = (1-x)^{\epsilon/2} F(x)$. When
$\epsilon = 0$, we find that $\Phi_0(x) = F(x)$ and $\Psi_0(x) =
-\frac{1}{s(s+1)} F(x) + {\rm const}$, such that $\Psi_0 \in X$ if
and only if $F(x) \in X$. The only set of eigenfunctions of the
Legendre equation (\ref{associated-Legendre-zero}) with $\epsilon =
0$ in $X$ is the set of Legendre polynomials $F = P_n(x)$ for $s =
n$ with $n \geq 0$ (see 8.91 on p. 973 in \cite{GR}). This set
corresponds to the eigenvalues (\ref{eigenvalues-1-zero}). Although
the zero eigenvalue $(s = n = 0)$ is excluded from the approach
above, it is still added to the spectrum by Lemma \ref{lemma-zero}.

When $\epsilon > 0$, the eigenfunction $\Psi_0(x)$ belongs to $X$
only if $F(x)$ has a regular behavior as $x = 1$ ($z = 0$). The
only solution of the hypergeometric equation
(\ref{hypergeometric-equation}) which is bounded as $x \to 1$ is
the hypergeometric function $F(z;\alpha,\beta,\gamma)$. (Indeed,
by 9.153 on p. 1001 of \cite{GR}, the other linearly independent
solution $F(x)$ has a singular behavior like $F(x) \sim
(1-x)^{-\epsilon/2}$ as $x \to 1$, which results in the divergence
$\Psi'_0(x) \sim (1-x)^{-1}$ as $x \to 1$, such that $\Psi_0
\notin X$.) By the identity 9.131 on p.998 of \cite{GR}, the
hypergeometric function $F(z;\alpha,\beta,\gamma)$ admits the
following behavior at the other singular point $x = -1$ ($z = 1$):
\begin{eqnarray}
\nonumber F(z;\alpha,\beta,\gamma) = \frac{\Gamma(\gamma)
\Gamma(\gamma-\alpha-\beta)}{\Gamma(\gamma-\alpha)
\Gamma(\gamma-\beta)}  F(1-z;\alpha,\beta,\alpha+\beta-\gamma+1)\\
\label{hypergeometric-expansion} + (1-z)^{\gamma-\alpha-\beta}
\frac{\Gamma(\gamma) \Gamma(\alpha+\beta-\gamma)}{\Gamma(\alpha)
\Gamma(\beta)}
F(1-z;\gamma-\alpha,\gamma-\beta,\gamma-\alpha-\beta+1),
\end{eqnarray}
where $\Gamma(z)$ is the Gamma function and
$$
z = \frac{1-x}{2}, \quad \alpha = \frac{\epsilon}{2} - s, \quad
\beta = \frac{\epsilon}{2} + s + 1, \quad \gamma =
\frac{\epsilon}{2} + 1.
$$
Since $\alpha + \beta - \gamma + 1 = \gamma$ and $\gamma - \alpha -
\beta + 1 = 1 - \frac{\epsilon}{2}$, the relation
(\ref{hypergeometric-expansion}) can be used only for $\epsilon < 2$
(the hypergeometric function $F(z;\alpha,\beta,\gamma)$ diverges for
$\gamma = -n$ with $n \geq 0$ integer).

It follows from (\ref{Psi-0-equation}) that the first term in
(\ref{hypergeometric-expansion}) leads the singular behavior of
$\Psi_0'(x) \sim (1+x)^{-1}$ as $x \to -1$ ($z \to 1$) if $\epsilon
\neq 0$, while the second term in (\ref{hypergeometric-expansion})
leads to the singular behavior $\Psi_0'(x) \sim (1+x)^{-\epsilon/2}$
as $x \to -1$ ($z \to 1$) if $\mu \neq 0$. Therefore, the
eigenfunction $\Psi_0(x)$ belongs to $X$ only if the first term in
(\ref{hypergeometric-expansion}) is removed which is only possible
if $\gamma - \alpha = 1 + s = -n$ or $\gamma - \beta = -s = -m$ with
integers $n,m \geq 0$. Both choices define the same set of
eigenvalues (\ref{eigenvalues-1-zero}) in the parametrization $\mu =
-s(s+1)$. Using another identity 9.131 on p.998 of \cite{GR},
\begin{equation}
\label{identity-hypergeometric} F(z;\alpha,\beta,\gamma) =
(1-z)^{\gamma-\alpha-\beta} F(z;\gamma-\alpha,\gamma-\beta,\gamma),
\end{equation}
we set $s = -1-n$ with $n \geq 1$, such that
$$
F\left(z;\frac{\epsilon}{2}+1+n,\frac{\epsilon}{2}-n,\frac{\epsilon}{2}
+1 \right) = (1-z)^{-\epsilon/2} F\left(z;-n,n+1, \frac{
\epsilon}{2}+1\right),
$$
where $F\left(z;-n,n+1,1 + \epsilon/2\right) \equiv \tilde{F}_n(x)$
is a polynomial of degree $n$, e.g.
\begin{eqnarray*}
\tilde{F}_0 = 1, \; \tilde{F}_1 = \frac{x + \sigma}{1 + \sigma}, \;
\tilde{F}_2 = \frac{3 x^2 + 3 \sigma x + \sigma^2
-1}{(1+\sigma)(2+\sigma)}, \; \tilde{F}_3 = \frac{15 x^3 + 15 \sigma
x^2 + (6 \sigma^2 - 9) x + \sigma (\sigma^2 -
4)}{(1+\sigma)(2+\sigma)(3+\sigma)},
\end{eqnarray*}
with $\sigma = \epsilon /2$. When $\epsilon = 0$ ($\sigma = 0$),
polynomials $\tilde{F}_n$ coincide with Legendre polynomials
$P_n(x)$ in 8.91 on p. 973 of \cite{GR}. The zero eigenvalue $(n =
0)$ is excluded from the construction but added to the spectrum by
Lemma \ref{lemma-zero}. When $\epsilon < 2$, the resulting
eigenfunction $\Psi_0(x)$ belongs to $X$. Since $\Psi_0'(x) \sim
(1+x)^{-\epsilon/2}$ as $x \to -1$, the resulting eigenfunction
$\Psi_0(x)$ does not belong to $X$ for $\epsilon \geq 2$.

We shall prove that no non-zero eigenvalues exist in $X$ for
$\epsilon \geq 2$. Using the identity
(\ref{identity-hypergeometric}), we transform the solution $F(x)$ to
the equivalent form $F(x) = (1+x)^{-\epsilon/2} \tilde{F}(x)$, where
$\tilde{F}(x)$ satisfies the hypergeometric equation
(\ref{hypergeometric-equation}) with new parameters
$$
z = \frac{1-x}{2}, \quad \tilde{\alpha} = \gamma - \alpha = 1 + s,
\quad \tilde{\beta} = \gamma - \beta = -s, \quad \tilde{\gamma} =
\gamma = \frac{\epsilon}{2} + 1.
$$
Up to a constant factor, $\tilde{F}(x)$ is represented by the
hypergeometric function $F(z;1+s,-s,1+\epsilon/2)$. It follows from
the ODE (\ref{Psi-0-equation}) that the eigenfunction $\Psi_0(x)$ is
related to $\tilde{F}(x)$ by
$$
-s(s+1) \Psi_0'(x) = \left( \frac{1 - x}{1+x} \right)^{\epsilon / 2}
\tilde{F}'(x).
$$
Since $\tilde{\alpha} + \tilde{\beta} - \tilde{\gamma} = -
\epsilon/2 < 0$ for $\epsilon > 0$, the hypergeometric series for
the function $F(z;1+s,-s,1+\epsilon/2)$ converges absolutely on
the entire interval $x \in [-1,1]$ ($z \in [0,1]$) (see 9.102 on
p.995 of \cite{GR}). Therefore, $\tilde{F}(x) \in C^2$ on $x \in
[-1,1]$ and $\tilde{F}'(-1)$ is well-defined. We shall prove that
$\tilde{F}'(-1) \neq 0$ for any $s \neq 0$ and $\epsilon \geq 2$.
It follows from the hypergeometric equation
(\ref{hypergeometric-equation}) with
$(\tilde{\alpha},\tilde{\beta},\tilde{\gamma})$ at $z = 1$ that
$$
s(s+1) \tilde{F}(-1) + \left( \frac{\epsilon}{2} - 1 \right)
\tilde{F}'(-1) = 0.
$$
If $\tilde{F}'(-1) = 0$, then $\tilde{F}(-1) = 0$ for any $s \neq
0$ and $\epsilon \geq 2$, and the only regular solution of the
hypergeometric equation (\ref{hypergeometric-equation}) is
$\tilde{F}(x) \equiv 0$. Therefore, $\tilde{F}'(-1) \neq 0$, and
therefore, $\Psi_0 \notin X$ for $\epsilon \geq 2$.
\end{proof}

\begin{proposition}
\label{proposition-epsilon-0-k-0} A complete spectrum of the
eigenvalue problem (\ref{zero-equation-stability}) with $0 < x_0 <
1$ and $\epsilon \geq 0$ in $X_0$ consists of simple isolated
eigenvalues $\mu$ with $\mu \in \mathbb{R}_-$.
\end{proposition}

\begin{proof}
When $\Psi_0 \in X_0$, the eigenfunction $\varphi(x)$ of the
associated Legendre equation (\ref{associated-Legendre-zero})
satisfies the Robin boundary conditions
$$
2(1-x_0^2) \varphi'(\pm x_0) + \epsilon \varphi(\pm x_0) = 0,
$$
while the eigenfunction $\chi(x)$ of the second-order ODE
(\ref{self-adjoint-form}) satisfies the Dirichlet boundary
conditions $\chi(\pm x_0) = 0$. Each eigenvalue problem is
self-adjoint with respect to these boundary conditions
\cite{Strauss}. Therefore, all eigenvalues $\mu$ of the regular
boundary-value problems are real-valued and isolated. Moreover,
these eigenvalues are negative due to the Green's identity
\cite{Strauss}:
\begin{equation}
\label{Green-identity-11} \mu \int_{-x_0}^{x_0} \varphi^2(x) dx =
- \int_{-x_0}^{x_0} (1-x^2) \left( \varphi'(x)\right)^2 dx -
\frac{\epsilon^2}{4} \int_{-x_0}^{x_0} \frac{\varphi^2(x)}{1-x^2}
dx < 0.
\end{equation}
These eigenvalues are also simple, since the Wronskian of any two
solutions with the Robin or Dirichlet boundary conditions is zero.
\end{proof}

\begin{proposition}
Let $\{ \mu_n \}_{n \geq 0}$ be isolated simple eigenvalues of the
eigenvalue problem (\ref{self-adjoint-form}) with $0 < x_0 < 1$
and Dirichlet boundary conditions $\chi(\pm x_0) = 0$ ordered as
$$
0 > \mu_0 > \mu_1 > ... > \mu_n > ...
$$
Then, $\lim\limits_{x_0 \to 1} \mu_n = - s_n(s_n+1)$, where
$$
s_n = 1 + n, \;\; \mbox{for} \;\; 0 \leq \epsilon \leq 2 \qquad
\mbox{and} \qquad s_n = \frac{\epsilon}{2} + n, \;\; \mbox{for} \;\;
\epsilon \geq 2
$$
with $n \geq 0$. \label{proposition-continuity-zero}
\end{proposition}

\begin{proof}
Singularity analysis of the second-order ODE
(\ref{self-adjoint-form}) shows that the solution $\chi(x)$ behaves
as
$$
\chi \to c_1^+ (1 - x)^{(\epsilon+2)/4} + c_2^+
(1-x)^{-(\epsilon+2)/4}, \quad \mbox{as} \;\; x \to 1
$$
and
$$
\chi \to c_1^- (1+x)^{(\epsilon-2)/4} + c_2^- (1+x)^{-(\epsilon
-2)/4}, \quad \mbox{as} \;\; x \to -1
$$
The ODE (\ref{self-adjoint-form}) admits a bounded (regular)
solution $\chi(x)$ on $x \in [-1,1]$ only if the singular
components are removed. This leads to the constraints $c_2^+ = 0$
and either $c_1^- = 0$ for $0 \leq \epsilon < 2$ or $c_2^- = 0$
for $\epsilon > 2$. It is explained in Proposition
\ref{proposition-x0-1-k-0} that the set $c_2^+ = 0$ and $c_1^- =
0$ for $0 \leq \epsilon < 2$ is equivalent to $s = m$ with $m \geq
0$, when the first term in the relation
(\ref{hypergeometric-expansion}) is removed and the hypergeometric
function $F(z;\tilde{\alpha},\tilde{\beta},\tilde{\gamma})$ is a
polynomial. Note that the zero eigenvalue $s = 0$ ($m = 0$) of the
problem (\ref{associated-Legendre-zero}) is excluded from the
spectrum of the problem (\ref{self-adjoint-form}), such that $s =
s_n = 1 + n$ with $n \geq 0$. On the other hand, the set $c_2^+ =
0$ and $c_2^- = 0$ for $\epsilon > 2$ is equivalent to $s = s_n =
\epsilon/2 + n$ with $n \geq 0$, when the second term in the
relation (\ref{hypergeometric-expansion}) is removed and the
hypergeometric function $F(z;\alpha,\beta,\gamma)$ is a
polynomial. Note that the first Darboux--Backlund transformation
(\ref{backlund-1}) implies that if $c_2^+ = 0$, then
$$
\varphi \to c_1^+ (1 - x)^{\epsilon/4}, \quad \mbox{as} \;\; x \to 1
$$
and
$$
\varphi \to c_1^- (1+x)^{\epsilon/4} + c_2^- (1+x)^{-\epsilon/4},
\quad \mbox{as} \;\; x \to -1.
$$
Recall that $\varphi(x) = (1-x^2)^{\epsilon/4} F(x)$. When $c_1^-
= 0$ ($0 \leq \epsilon < 2$), $F(x)$ is singular like $F(x) \to
(1+x)^{-\epsilon/2}$ as $x \to -1$ in accordance with the relation
(\ref{identity-hypergeometric}). When $c_2^- = 0$ ($\epsilon >
2$), $F(x)$ is bounded as $x \to -1$. The marginal case $\epsilon
= 2$ corresponds to the case when $\chi(x)$ has a bounded and
logarithmically  growing components as $x \to -1$. The logarithmic
growth is excluded if $s_n = 1 + n$ with $n \geq 1$, which is the
border between the two spectra at $\epsilon = 2$.

When $s = s_n$ and $\epsilon \neq 2$, the eigenfunction $\chi(x)$
of the formally self-adjoint problem (\ref{self-adjoint-form})
satisfies the Dirichlet boundary conditions $\lim_{x \to \pm 1}
\chi(x) = 0$. When $\epsilon = 2$, the eigenfunction $\chi(x)$ is
bounded at $x = -1$ and zero at $x = 1$. In either case,
convergence and uniqueness of continuations from eigenvalues of
the regular Dirichlet problem (\ref{self-adjoint-form}) with $x_0
< 1$ to eigenvalues of the singular boundary-value problem
(\ref{self-adjoint-form}) with $x_0 = 1$ is proved by Theorem 5.3
of \cite{zettle}.
\end{proof}

\begin{remark}
{\rm Bounded (for $\epsilon = 2$) and decaying (for $\epsilon
> 2$) eigenfunctions $\chi(x)$ of the self-adjoint problem
(\ref{self-adjoint-form}) with $x_0 = 1$ for eigenvalues $\mu =
-s_n(s_n+1)$ with $s_n = \epsilon/2 + n$ violate the conditions
(\ref{conditions-chi-X}). Indeed, one can check that the limit in
(\ref{conditions-chi-X}) as $x \to -1$ is non-zero (proportional
to $c_1^-$) and the integral in (\ref{conditions-chi-X}) hence
diverges. Therefore, the eigenvalues of the self-adjoint problem
(\ref{self-adjoint-form}) for $\epsilon \geq 2$ do not correspond
to eigenvalues of the original problem
(\ref{zero-equation-stability}) in space $\Psi_0 \in X$, in
agreement with Proposition \ref{proposition-x0-1-k-0}. We also
note that if one consider a generalized conditions for $\chi(x)$
with
\begin{equation}
\left| \lim_{x \to \pm 1} \left( \frac{1-x}{1+x}
\right)^{\epsilon/4} \sqrt{1-x^2} \chi(x) \right| < \infty,
\end{equation}
the spectrum of the self-adjoint problem (\ref{self-adjoint-form})
is not defined since $c_2^+ \neq 0$ and the eigenvalue problem is
not complete. }
\end{remark}

\begin{remark}
{\rm Theorem \ref{theorem-2} covers the entire parameter domain $0
< x_0 \leq 1$ and $\epsilon \geq 0$. However, there is an
interesting problem with convergence of eigenvalues of the
associated Legendre equation (\ref{associated-Legendre-zero}) in
the limit $x_0 \to 1$. While the eigenvalues with $0 < x_0 < 1$
are expected to converge to the eigenvalues in
(\ref{eigenvalues-1-zero}) for $0 \leq \epsilon < 2$, no
eigenvalues with the eigenfunctions $\Psi \in X$ exist for
$\epsilon \geq 2$. Convergence of eigenvalues of the linearized
problem (\ref{zero-equation-stability}) as $x_0 \to 1$ will be
computed in Section 6 numerically. }
\end{remark}

\section{Numerical computations of eigenvalues for $k \neq 0$}

In order to illustrate distribution of eigenvalues in Propositions
\ref{proposition-x0-1}, \ref{proposition-epsilon-0} and
\ref{proposition-continuity} and to investigate eigenvalues in the
domain $0 < x_0 < 1$ and $\epsilon > 0$ in Remark \ref{remark-k},
we develop a numerical method based on power series expansions.
Since $x = 0$ is an ordinary point and $x = \pm 1$ are regular
singular points of the system (\ref{stability-problem}), the power
series expansions of the functions $\Psi_k(x)$ and $\Phi_k(x)$ in
powers of $x$ converge uniformly and absolutely for $|x| < 1$. The
numerical method is based on truncation of the power series.

Let $\mu \in \mathbb{C}$ be parameterized by $\mu = -s(s+1)$, $s
\in \mathbb{C}$. Due to the symmetry, it is sufficient to consider
the domain $\{ s \in \mathbb{C} : \; {\rm Re}(s) \geq -\frac{1}{2}
\}$. The stability domain ${\rm Re}(\mu) < 0$ corresponds to the
domain
\begin{equation}
\label{stability-boundary} \left\{ s \in \mathbb{C} : \quad |{\rm
Im}(s)| < \sqrt{{\rm Re}(s)( {\rm Re}(s) + 1)}, \;\; {\rm Re}(s) >
0 \right\}.
\end{equation}
Consider the power series with separated even and odd terms:
\begin{eqnarray}
\label{Psi-series-k} \Psi_k(x) & = & \sum_{m \geq 0} c_m x^{2m} +
\sum_{m \geq 0} d_m x^{2m+1}, \\
\label{Phi-series-k} \Phi_k(x) & = & \sum_{m \geq 0} a_m x^{2m} +
\sum_{m \geq 0} b_m x^{2m+1},
\end{eqnarray}
where the starting coefficients $(a_0,b_0,c_0,d_0)$ are parameters.
Substituting (\ref{Phi-series-k}) into the second equation of the
system (\ref{stability-problem}) we find that $(a_1,b_1)$ are
defined separately as
\begin{eqnarray}
\label{recurrence1-app} a_1 & = & \frac{(k^2 - s(1+s)) a_0 -
\epsilon b_0}{2}, \\
\label{recurrence2-app} b_1 & = & \frac{(k^2 + 2 - s(s+1)) b_0 -
\epsilon 2 a_1}{6},
\end{eqnarray}
while the coefficients $\{ a_m,b_m \}_{m \geq 2}$  are defined
uniquely from the recurrence equations:
\begin{eqnarray}
\nonumber a_{m+2} & = & \frac{(k^2 - s(s+1) + 2(2m+2)^2) a_{m+1} +
(s(s+1)-2m (2m+1)) a_m}{(2m+4)(2m+3)} \\
\label{recurrence1-app-2} & \phantom{t} & \phantom{texttext} \frac{-
\epsilon (2m+3) b_{m+1} + \epsilon (2m+1) b_m}{(2m+4)(2m+3)}, \\
\nonumber b_{m+2} & = & \frac{(k^2 - s(s+1) + 2 (2m+3)^2) b_{m+1} +
(s(s+1)- (2m+2) (2m+1)) b_m}{(2m+5)(2m+4)} \\
\label{recurrence2-app-2} & \phantom{t} & \phantom{texttext}
\frac{-\epsilon (2m+4) a_{m+2} + \epsilon (2m+2)
a_{m+1}}{(2m+5)(2m+4)}.
\end{eqnarray}
We note that the initial equations
(\ref{recurrence1-app})--(\ref{recurrence2-app}) follow from the
recurrence equations
(\ref{recurrence1-app-2})--(\ref{recurrence2-app-2}) for $m = -1$
with $a_{-1} = b_{-1} = 0$.

Substituting (\ref{Psi-series-k}) into the first equation of the
system (\ref{stability-problem}) we find that the coefficients $\{
c_m,d_m \}_{m \geq 2}$  are defined from the coefficients $\{
a_m,b_m \}_{m \geq 0}$ by the recurrence equations:
\begin{eqnarray}
\label{recurrence3-app-2} c_{m+2} & = & \frac{(k^2 + 2(2m+2)^2)
c_{m+1} - 2m (2m+1) c_m + a_{m+1} - a_m}{(2m+4)(2m+3)} \\
\label{recurrence4-app-2} d_{m+2} & = & \frac{(k^2 + 2 (2m+3)^2)
d_{m+1} -(2m+2)(2m+1) d_m + b_{m+1} - b_m}{(2m+5)(2m+4)}.
\end{eqnarray}
The initial equations for $(c_1,d_1)$ follow from the recurrence
equations (\ref{recurrence3-app-2})--(\ref{recurrence4-app-2}) for
$m = -1$ with $a_{-1} = b_{-1} = c_{-1} = d_{-1} = 0$.

The boundary conditions in (\ref{bc1}) lead to the equations
\begin{eqnarray} \label{boundary-condition-1-k} \sum_{m \geq 0}
c_m x_0^{2m} = 0, \quad \sum_{m \geq 0} d_m x_0^{2m} = 0, \quad
\sum_{m \geq 0} (2m) c_m x_0^{2m} = 0, \quad \sum_{m \geq 0}
(2m+1) d_m x_0^{2m} = 0.
\end{eqnarray}
There exists a linear map from $(a_0,b_0,c_0,d_0) \in
\mathbb{C}^4$ parametrized by $s \in \mathbb{C}$ to the sequence
$\{ a_m,b_m,c_m,d_m \}_{m \in \mathbb{N}}$. Therefore, the
boundary conditions (\ref{boundary-condition-1-k}) are equivalent
to the homogeneous system $A_k(s) {\bf x} = {\bf 0}$, where ${\bf
x} = (a_0,b_0,c_0,d_0)^T \in \mathbb{C}^4$ and $A_k(s)$ is a
$4$-by-$4$ matrix computed from the entries of
(\ref{boundary-condition-1-k}). The matrix $A_k(s)$ depends on $s
\in \mathbb{C}$ and $k \in \mathbb{N}$, as well as parameters
$x_0$ and $\epsilon$. If the power series are truncated at the
$M$-th term, the matrix $A_k(s)$ depends also on $M$. Eigenvalues
$\mu = -s(s+1)$ of the system (\ref{stability-problem}) in
(\ref{bc1}) are {\em equivalent} to roots $s$ of the determinant
equation
\begin{equation}
F_k(s;x_0,\epsilon,M) = {\rm det}(A_k(s)).
\end{equation}

\begin{figure}[htbp]
\begin{center}
\includegraphics[height=8cm]{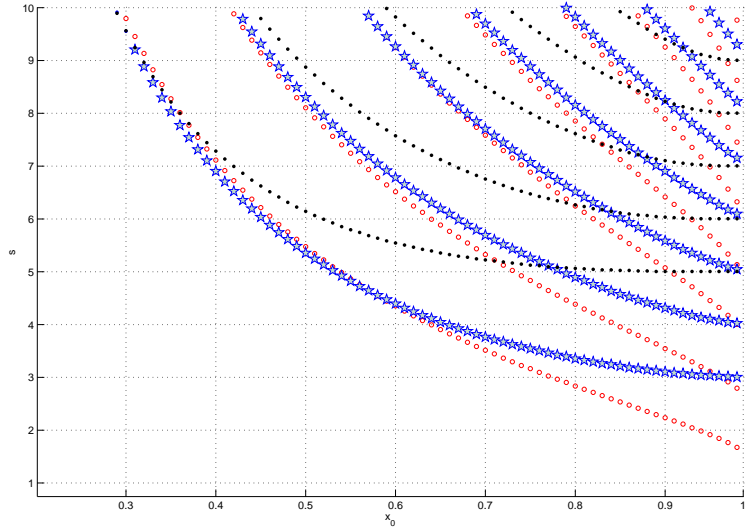}
\end{center}
\caption{First few roots $s$ of $F_k(s;x_0,\epsilon,M)$ versus
$x_0$ for $\epsilon = 0$ and $M = 150$: $k = 1$ (circles), $k = 3$
(stars) and $k = 5$ (dots).} \label{fig-s-x0a}
\end{figure}

Numerical results of computations of roots of the function
$F_k(s;x_0,\epsilon,M)$  are shown on Figures 1--5. Figure
\ref{fig-s-x0a} show first few roots $s$ of
$F_k(s;x_0,\epsilon,M)$ with $k = 1, 3, 5$ versus $x_0$ for
$\epsilon = 0$ and $M = 150$. In agreement with Proposition
\ref{proposition-continuity}, the roots converge as $x \to 1$ to
the values $s_n = \sigma + n$ with $\sigma = |k|$ and $n \geq 0$.
We can see that the convergence is excellent for $k = 3$ and $k =
5$ but it is worse for $k = 1$ in the sense that the roots at $x_0
= 0.99$ are still far from the values $s_n$. This feature is
explained by the decay of the eigenfunctions $(\Phi_k,\Psi_k)$ of
the system (\ref{stability-problem}) on $x \in [-1,1]$. Indeed, it
follows from Proposition \ref{proposition-x0-1} that $\Phi_k \sim
(1 - x^2)^{\sigma/2}$ and $\Psi_k \sim (1-x^2)^{1 + \sigma/2}$ as
$x \to \pm 1$ for $\epsilon = 0$ and $|k| \geq 1$. Therefore, the
derivative of $\Phi_k(x)$ is bounded as $x \to \pm 1$ for $|k|
\geq 2$ and unbounded for $|k| = 1$. In the latter case,  the
power series expansions (\ref{Psi-series-k})--(\ref{Phi-series-k})
diverge in the limit $x_0 \to 1$ and the numerical approximation
is not accurate for $x_0$ close to $1$.

\begin{figure}[htbp]
\begin{center}
\includegraphics[height=8cm]{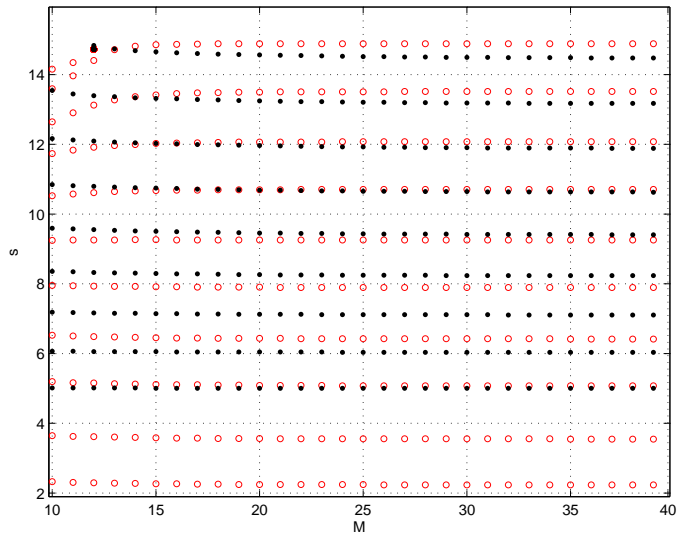}
\end{center}
\caption{Convergence of roots $s$ versus $M$ for $\epsilon = 0$
and $x_0 = 0.9$: $k = 1$ (circles) and $k = 5$ (dots). }
\label{fig-s-M}
\end{figure}

Figure \ref{fig-s-M} shows first few roots $s$ with $k = 1,5$
versus $M$ for $\epsilon = 0$ and $x_0 = 0.9$. We can see that the
roots quickly converge to constant values, which are taken as
approximations of real roots when $M = 150$ in the remainder of
the figures. The numerical error for large values of $M$ consists
of three sources: truncation of the power series, root finding
algorithms, and rounding entries of the matrix $A_k(s)$ when a
number $x_0$ with $x_0 < 1$ is evaluated at a large power $x_0^M$.
While the first two sources can be reduced to any desired degree,
the last source represents an irremovable obstacle on getting
accurate approximations when $M$ gets large.

\begin{figure}[htbp]
\begin{center}
\includegraphics[height=8cm]{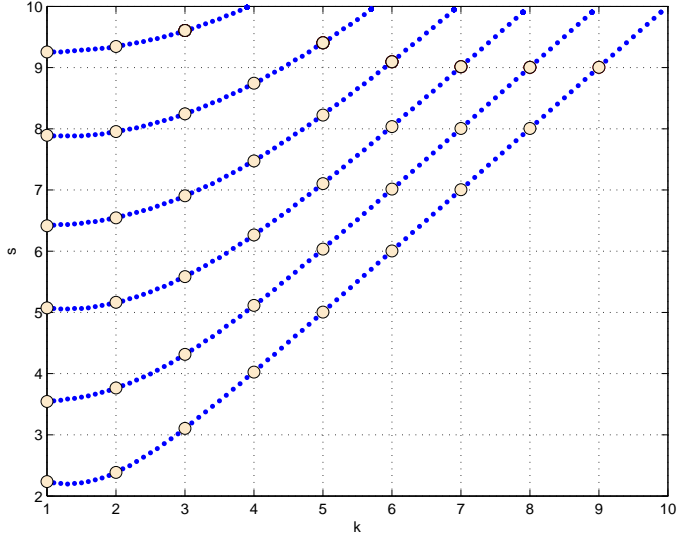}
\end{center}
\caption{First six roots $s$ versus $k$ for $\epsilon = 0$, $x_0 =
0.9$, and $M = 150$. } \label{fig-s-k}
\end{figure}

Figure \ref{fig-s-k} shows the first six roots $s$ versus $k$ for
$\epsilon = 0$, $x_0 = 0.9$, and $M = 150$. We observe two
properties from this figure: the values of $s$ becomes larger for
larger values of $k$ (e.g. the eigenvalues $\mu$ becomes more and
more negative) and the roots $s$ approach to the integer values
for larger values of $k$ even when $x_0 = 0.9$ is not close to
$x_0 = 1$.

\begin{figure}[htbp]
\begin{center}
\includegraphics[height=8cm]{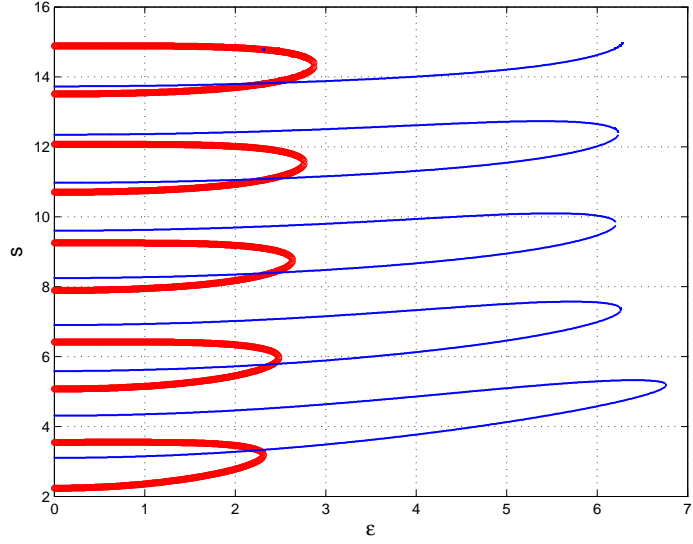}
\end{center}
\caption{First few roots $s$ versus $\epsilon$ for $x_0 = 0.9$ and
$M = 150$: $k = 1$ (bolded curve) and $k = 3$ (thin curve).}
\label{fig-s-eps}
\end{figure}

Figure \ref{fig-s-eps} show the first few roots $s$ with $k = 1,3$
versus $\epsilon$ for $x_0 = 0.9$ and $M = 150$. Although the
roots are real for small values of $\epsilon$ in agreement to
Proposition \ref{proposition-epsilon-0}, they coalesce for larger
values of $\epsilon$. After two roots merge, they split into
complex domain and complex values of $s$ are not shown on Figure
\ref{fig-s-eps}. It is seen from this figure that the roots with
larger values of $k$ coalesce for larger values of $\epsilon$.

\begin{figure}[htbp]
\begin{center}
\includegraphics[height=6cm]{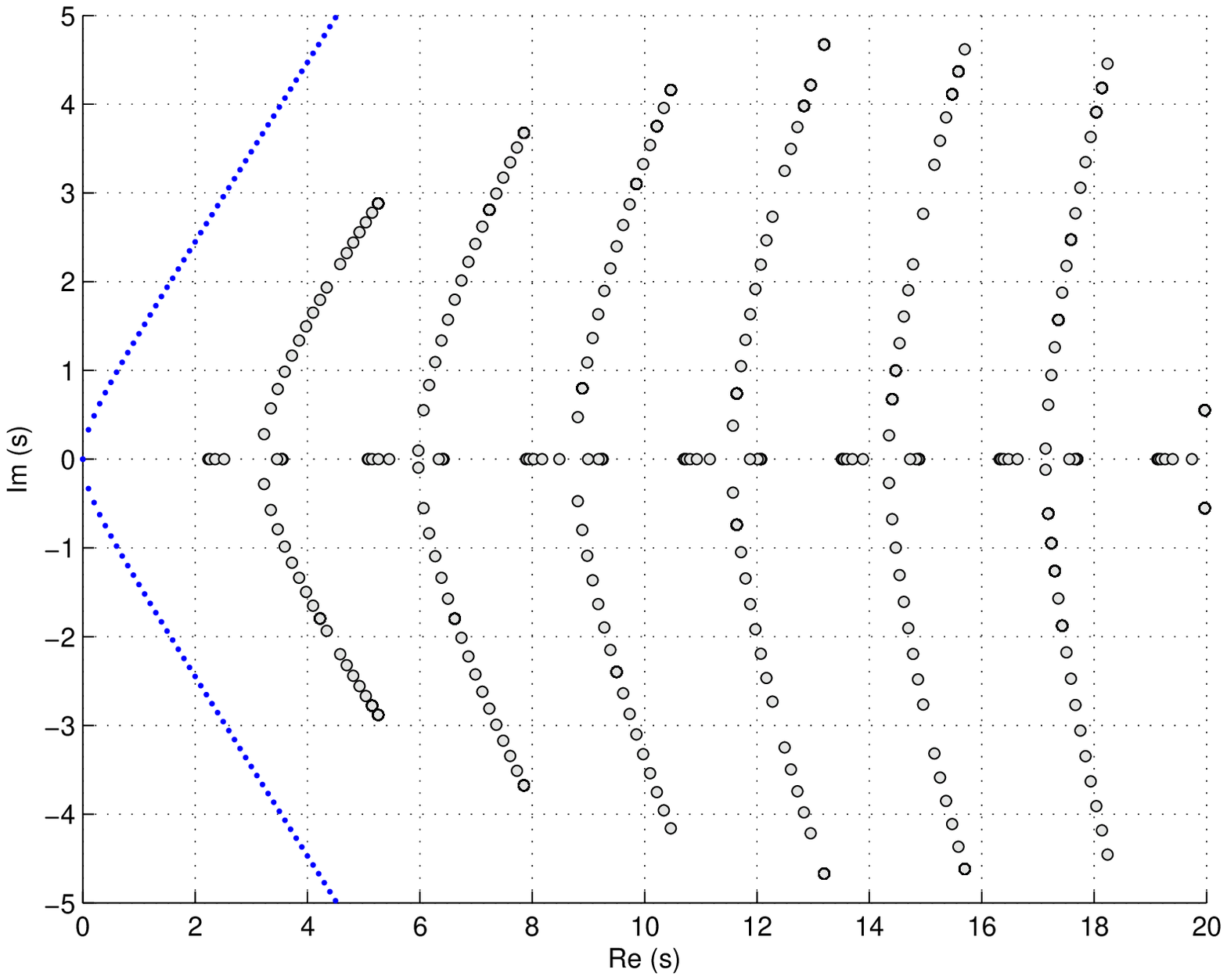}
\includegraphics[height=6cm]{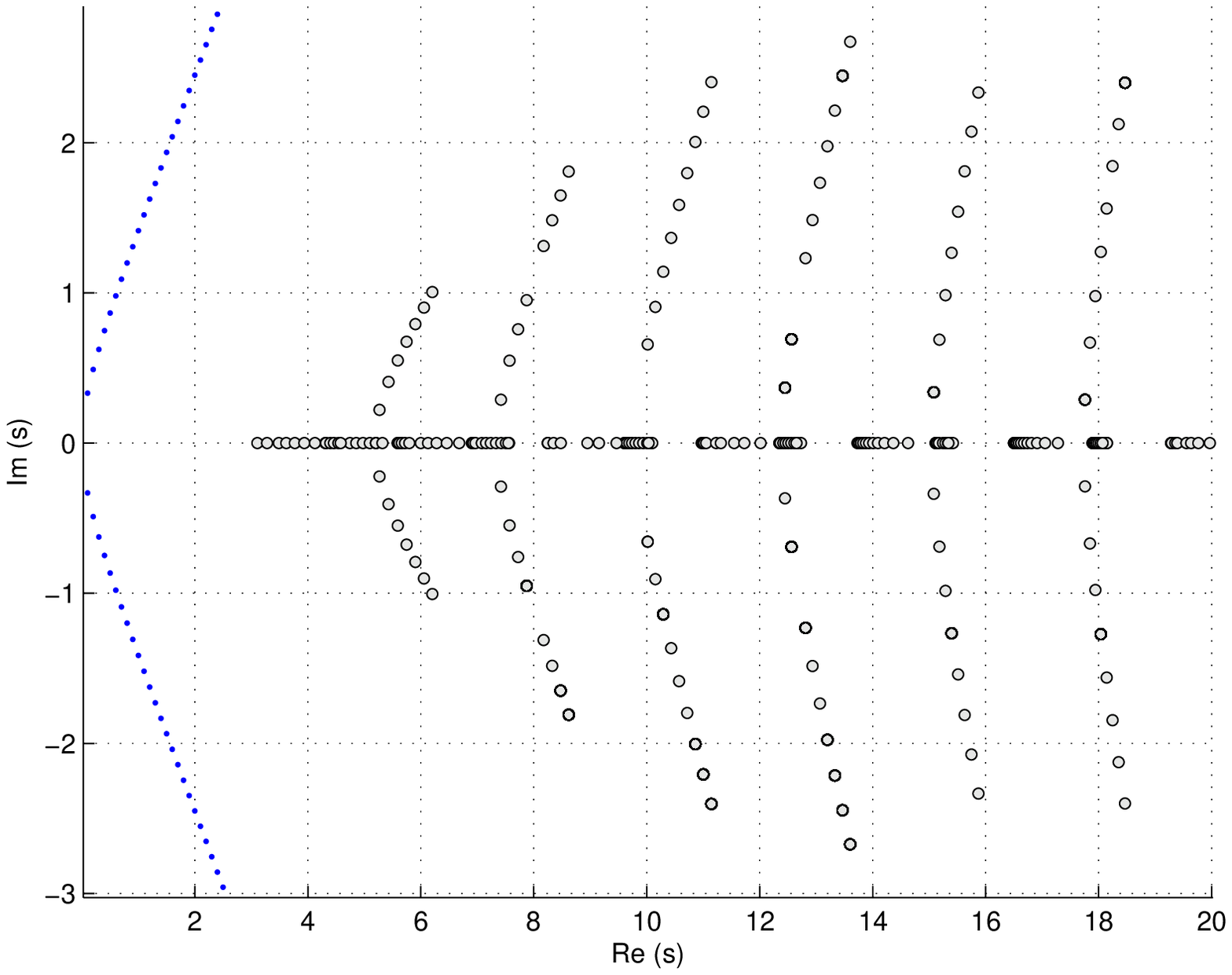}
\end{center}
\caption{Complex roots $s$ for $k = 1$ (left) and $k = 3$ (right),
$x_0 = 0.9$ and $M = 150$ when parameter $\epsilon$ transverses in
the interval $0 \leq \epsilon \leq 12$. The dotted curve shows the
boundary of the stability domain (\ref{stability-boundary}).}
\label{fig-eps-compl}
\end{figure}

Figure \ref{fig-eps-compl} shows the spectrum of complex roots $s$
with $k = 1,3$ for $x_0 = 0.9$, $M = 150$, and different values of
$0 \leq \epsilon \leq 12$. The boundary of the stability domain
(\ref{stability-boundary}) is shown by the dotted curve. We can
see that roots $s$ remain in the stability domain after they
bifurcate off the real axes.

\section{Numerical computations of eigenvalues for $k = 0$}

We approximate eigenvalues of the system
(\ref{zero-equation-stability}) with power series solutions
explained in Section 5. The solution for $\Psi_0(x)$ and
$\Phi_0(x)$ is represented by the power series
(\ref{Psi-series-k})--(\ref{Phi-series-k}), where the starting
coefficients $(a_0,b_0,c_0,d_0)$ are parameters, while the
coefficients $\{ a_m,b_m,c_m,d_m \}_{m \in \mathbb{N}}$ are
defined uniquely from the recurrence equations. It follows from
the ODE (\ref{stability-zero-problem}) that the set $\{ a_m,b_m
\}_{m \in \mathbb{N}}$ is uncoupled from the other coefficients
but it is defined by the unknown value of the parameter $s$:
\begin{eqnarray}
\label{recurrence1} a_{m+1} & = & \frac{(2m-s)(2m+1+s) a_m -
\epsilon (2m+1) b_m}{(2m+2)(2m+1)}, \\
\label{recurrence2} b_{m+1} & = & \frac{(2m+1-s)(2m+2+s) b_m -
\epsilon (2m+2) a_{m+1}}{(2m+3)(2m+2)}.
\end{eqnarray}
Given $(a_0,b_0)$ and the value for $s$, the recurrence equation
(\ref{recurrence1}) gives the value of $a_1$ and then the recurrence
equation (\ref{recurrence2}) defines the value of $b_1$, and so on.
It follows from the first equation of the system
(\ref{zero-equation-stability}) that the set $\{ c_m,d_m \}_{m \in
\mathbb{N}}$ is defined by the set $\{ a_m,b_m \}_{m \in
\mathbb{N}}$ but it is independent of $s$:
\begin{eqnarray}
\label{recurrence3} c_{m+1} & = & \frac{(2m)(2m+1) c_m + a_m}{(2m+2)(2m+1)}, \\
\label{recurrence4} d_{m+1} & = & \frac{(2m+1)(2m+2) d_m +
b_m}{(2m+3)(2m+2)}.
\end{eqnarray}
Finally, it follows from the second equation of the system
(\ref{zero-equation-stability}) that there exist two initial
equations:
\begin{eqnarray*}
b_0 + \epsilon a_0 & = & -s(s+1) d_0, \\
2 a_1 + \epsilon b_0 & = & -2 s(s+1) c_1
\end{eqnarray*}
in addition to the system (\ref{recurrence1})--(\ref{recurrence2}).
When $s \neq 0$, we can solve the initial equations as
$$
b_0 = -\epsilon a_0 - s(s+1) d_0, \qquad c_1 = \frac{a_0}{2},
$$
such that the only independent parameters are $(a_0,d_0)$. We also
note that the parameter $c_0$ is trivial since $\Psi_0(x)$ is
defined up to the addition of an arbitrary constant.

The boundary conditions in (\ref{bc1-zero}) lead to the equations:
\begin{eqnarray}
\label{boundary-condition-1} \sum_{m \geq 0} (2m) c_m x_0^{2m} = 0,
\qquad \sum_{m \geq 0} (2m+1) d_m x_0^{2m} = 0.
\end{eqnarray}
There exists a linear map from $(a_0,d_0) \in \mathbb{C}^2$
parameterized by $s \in \mathbb{C}$ to the sequence $\{
a_m,b_m,c_m,d_m \}_{m \in \mathbb{N}}$. Therefore, the boundary
conditions (\ref{boundary-condition-1}) are equivalent to the
homogeneous system $A_0(s) {\bf x} = {\bf 0}$, where ${\bf x} =
(a_0,d_0)^T \in \mathbb{C}^2$ and $A_0(s)$ is a $2$-by-$2$ matrix
which depends on $s \in \mathbb{C}$, parameters $x_0$ and
$\epsilon$, and integer $M$ for truncation of power series.
Eigenvalues $\mu = -s(s+1)$ of the system
(\ref{zero-equation-stability}) in (\ref{bc1-zero}) are {\em
equivalent} to roots $s$ of the determinant equation
\begin{equation}
F_0(s;x_0,\epsilon,M) = {\rm det}(A_0(s)).
\end{equation}

\begin{figure}[htbp]
\begin{center}
\includegraphics[height=8cm]{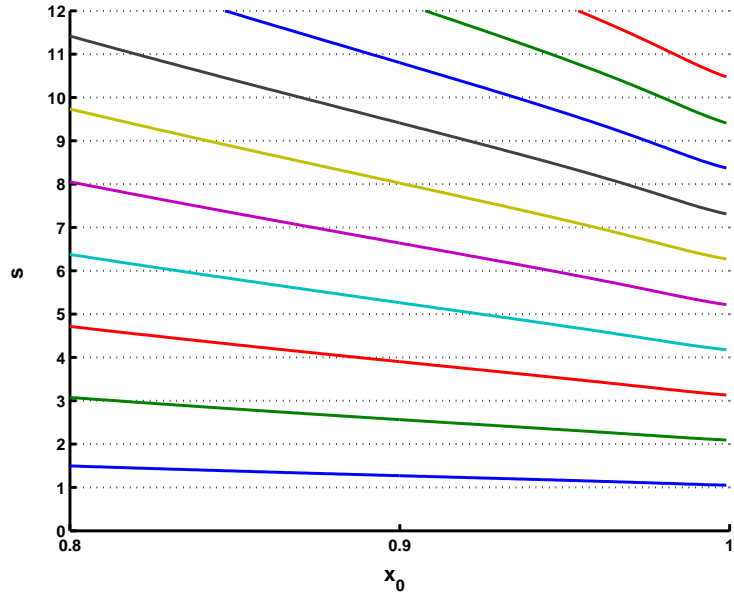}
\end{center}
\caption{First ten eigenvalues of the problem
(\ref{zero-equation-stability}) for $\epsilon = 1$ and $M = 100$.}
\label{fig-k0-1}
\end{figure}

Figure \ref{fig-k0-1} represents the first ten eigenvalues $s$
versus $x_0$ for $\epsilon = 1$ and $M = 100$. In agreement with
Proposition \ref{proposition-x0-1-k-0}, the roots converge to the
integer values in the limit $x_0 \to 1$. Since the convergence of
power series becomes slower with $M$ for $x_0 \neq 1$, there is a
gap between the last numerical data and the value $x_0 = 1$. We
also note that the numerical accuracy of the limiting eigenvalues
(\ref{eigenvalues-1-zero}) becomes worse for larger eigenvalues.

\begin{figure}[htbp]
\begin{center}
\includegraphics[height=8cm]{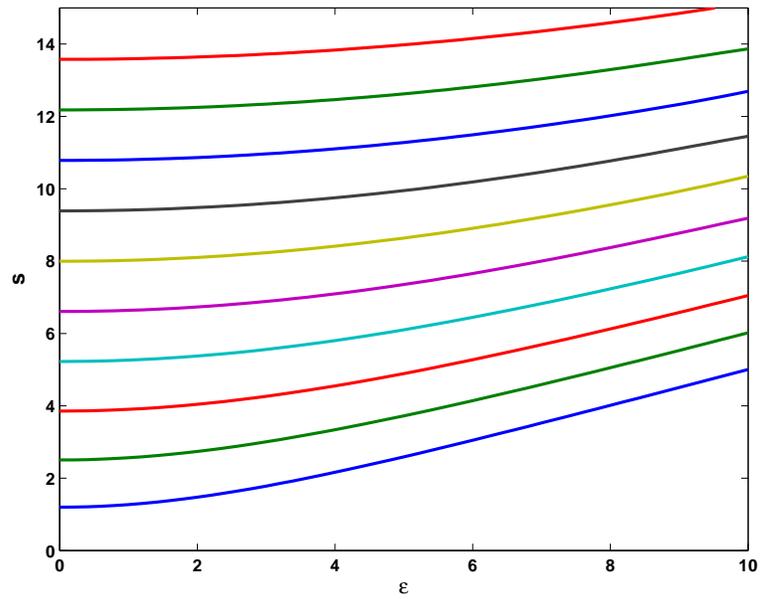}
\end{center}
\caption{First ten eigenvalues of the problem
(\ref{zero-equation-stability}) for $x_0 = 0.9$ and $M = 100$.}
\label{fig-k0-2}
\end{figure}

Figure \ref{fig-k0-2} represents the first ten eigenvalues $s$
versus $\epsilon$ for $x_0 = 0.9$ and $M = 100$. It is obvious
that the eigenvalues remain real in agreement with Proposition
\ref{proposition-epsilon-0-k-0}.

\begin{figure}[htbp]
\begin{center}
\includegraphics[height=8cm]{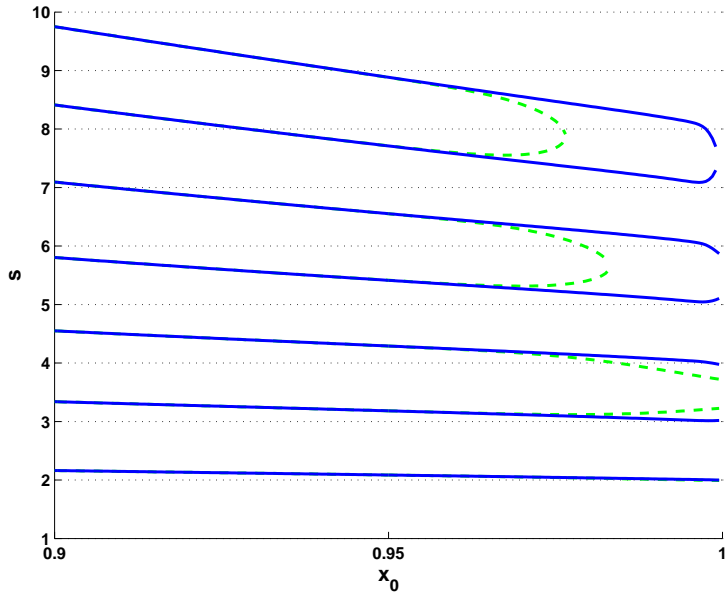}
\end{center}
\caption{Convergence of eigenvalues of the problem
(\ref{zero-equation-stability}) for $\epsilon = 4$ and two values
of $M = 100$ (dashed curve) and $M = 1000$ (solid curves). }
\label{fig-k0-3}
\end{figure}

Figure \ref{fig-k0-3} represents the first seven eigenvalues $s$
versus $x_0$ for $\epsilon = 4$ and two values of $M = 100$
(dashed curves) and $M = 1000$ (solid curves). In agreement with
Proposition \ref{proposition-continuity-zero}, the roots converge
to their limiting values which are not eigenvalues of the problem
(\ref{zero-equation-stability}) in space (\ref{bc2-zero}). We also
note limitations of the numerical methods based on truncations of
the power series. True limits can only be recovered if too many
terms of the power series are taken into accounts which leads to
long computational time and large round-off errors of numerical
computations. The effects of slow convergence and truncations of
power series lead to coalescence of real eigenvalues and their
splitting to the complex plane, which is not observed if the
values of $M$ are large enough.

\section{Discussions}

We have shown analytically that the stationary flow on the sphere
is asymptotically stable whatever the Reynolds number may occur.
This result is relevant for the flow of a viscous fluid (e.g. oil)
over a sphere (e.g. a metal ball). We have also found that the
linearized operator for symmetry-preserving perturbations has void
spectrum in the energy space for sufficiently large Reynolds
numbers. One can show by direct analysis that the full system
(\ref{2.1})--(\ref{2.3}) reduces to a scalar linear equation for
symmetry-preserving ($\phi$-independent) solutions:
\begin{equation}
\label{linear-equation-longitudinal} \frac{\partial
v_{\phi}}{\partial t} + \frac{1}{\sin \theta} \Delta_0 v_{\phi} =
\nu \frac{\partial}{\partial \theta} \Delta_0 v_{\phi },
\end{equation}
where $\Delta_0$ is given by (\ref{Laplacian-k}) for $k = 0$. When
$v_{\theta}(\theta,t) = - \Psi_0'(\theta) e^{\lambda t}$, the
linear equation (\ref{linear-equation-longitudinal}) reduces to
the linear eigenvalue problem (\ref{zero-equation}) which has no
eigenvalues in the space of square integrable functions
$\int_0^{\pi} \left( \Psi_0'(\theta)\right)^2 \sin \theta d \theta
< \infty$ when $\nu \leq \frac{1}{2}$ ($\epsilon \geq 2$).
Implications of this result to the well-posedness of the Cauchy
problem for the linear time-dependent equation
(\ref{linear-equation-longitudinal}) with $\nu \leq \frac{1}{2}$
remain unclear.

We have also shown analytically and numerically that the
stationary flow on the truncated spherical layer is asymptotically
stable and all isolated eigenvalues are real for small Reynolds
numbers and complex for large Reynolds numbers. The eigenvalues
are always real for symmetry-preserving perturbations. The
truncated spherical layer can be used to model the ice melting in
Arctics due to global warming, when the near-stationary flow of
ocean water moves from Arctics to Antarctica. We note however that
the model of two-dimensional Navier--Stokes equations on sphere
considered in this paper does not include the Earth's rotation,
the gravity force, and the location of continents.

{\bf Acknowledgement.} The authors thank Marina Chugunova and
Bartosz Protas for useful discussions and remarks. The work was
supported by the PREA and NSERC Discovery grants.

\end{document}